\documentclass [PhD] {uclathes}
\usepackage{amsmath}
\usepackage{amsfonts} 
\usepackage{amsthm} 
\usepackage{amssymb}
\usepackage{multirow}
\usepackage{mathtools}
\usepackage{hyperref}
\usepackage{tikz-cd}
\usepackage{graphicx}
\usepackage{array}
\usepackage[all,cmtip]{xy}

\newtheorem{theorem}{Theorem}[section]
\newtheorem{corollary}[theorem]{Corollary}
\newtheorem{proposition}[theorem]{Proposition}
\newtheorem{lemma}[theorem]{Lemma}
\newtheorem{construction}[theorem]{Construction}

\theoremstyle{definition}
\newtheorem{definition}[theorem]{Definition}

\theoremstyle{remark}

\newtheoremstyle{TheoremNum}
{7pt}{7pt}              %%% space between body and thm
{\itshape}                      %%% Thm body font
{}                              %%% Indent amount (empty = no indent)
{\bfseries}                     %%% Thm head font
{.}                             %%% Punctuation after thm head
{ }                             %%% Space after thm head
{\thmname{#1}\thmnote{ \bfseries #3}}%%% Thm head spec
\theoremstyle{TheoremNum}
\title          {\textnormal{\normalsize{The mixed Tate property of reductive groups}}}
\author         {\textnormal{\normalsize{Yehonatan Sella}}}
\department     {Mathematics}
\degreeyear     {2017}

\chair          {Burt Totaro}
\member         {Raphael Alexis Rouquier}
\member         {Alexander Sergee Merkurjev}
\member         {Paul Balmer}

\dedication{To my parents,\\ who have encouraged me to pursue my imagination}     

\acknowledgments {I am grateful to Burt Totaro for his help, guidance and encouragement as my advisor. I also thank my colleagues in the UCLA math department, especially Laure Flapan, Fei Xie and Bon-Soon Lin, who have been with me throughout the learning process and made it all the more enjoyable. Finally, I am grateful to Renee Bell who has been by my side throughout my PhD, even from afar.

The work in this thesis was partly supported by a Beckenbach Fellowship from UCLA.}

\vitaitem   {2012}
                {B.A.~(Mathematics),
                University of California, Berkeley.}

\abstract{This thesis is concerned with the mixed Tate property of reductive algebraic groups $G$, which in particular guarantees a Chow Kunneth property for the classifying space $BG$. Toward this goal, we first refine the construction of the compactly supported motive of a quotient stack. 

In the first section, we construct the compactly supported motive $M^c(X)$ of an algebraic space $X$ and demonstrate that it satisfies expected properties, following closely Voevodsky's work in the case of schemes. 

In the second section, we construct a functorial version of Totaro's definition of the compactly supported motive $M^c([X/G])$ for any quotient stack $[X/G]$ where $X$ is an algebraic space and $G$ is an affine group scheme acting on it. A consequence of functoriality is a localization triangle for these motives. 

In the third section, we study the mixed Tate property for the classical groups as well as the exceptional group $G_2$. For these groups, we demonstrate that all split forms satisfy the mixed Tate property, while exhibiting non-split forms that do not. Finally, we prove that for any affine group scheme $G$ and normal split unipotent subgroup $J$ of $G$, the motives $M^c(BG)$ and $M^c(B(G/J))$ are isomorphic.}

\begin {document}
\makeintropages

%%%%%%%%%%%%%%%%%%%%%%%%%%%%%%%%%%%%%%%%%%%%%%%%%%%%%%%%%%%%%%%%%%%%%%
%
% Ordinarily each chapter (at least) is in a separate file.
%

\chapter{Introduction}
This paper studies the Chow Kunneth property of classifying spaces of algebraic groups.

Given a topological group $G$, its classifying space is a space $BG$ together with a principal $G$-bundle $\pi:EG\to BG$ such that any principal $G$-bundle over a $CW$-complex $X$ is a pullback of the bundle $\pi$ by a morphism $X\to BG$ which is unique up to homotopy.

In the context of algebraic geometry, given an affine group scheme $G$ over a field $k$, the analogous concept is the quotient stack $[Spec k/G]$. By definition, for any $k$-scheme $X$, the groupoid of morphisms $X\to [Spec k/G]$ is the groupoid of all $G$-bundles over $X$.

Burt Totaro defined the Chow groups of classifying spaces, $CH^i(BG)$, by approximating $BG$ up to homotopy by a scheme which agrees with it modulo codimension $i$. Namely, one can find a representation $V$ of $G$ and a $G$-stable closed subset $S\subset V$  of codimension larger than $i$ such that $(V-S)/G$ is a scheme. Totaro defined $CH^i(BG)=CH^i((V-S)/G)$\cite{burtchow}. This definition was generalized to arbitrary quotient stacks $[X/G]$ by Edidin-Graham\cite{edidingraham}, where $X$ is an algebraic space and $G$ an affine group scheme acting on it. The Chow groups $CH^i([X/G])$ may be thought of as the $G$-equivariant Chow groups of $X$.

Classifying spaces provide a fruitful setting for the calculation of Chow groups, as they combine geometric with group- and representation-theoretic information, and are functorial with respect to groups.

One striking property of classifying spaces is that, while a $k$-scheme $X$ often does not satisfy the Chow-Kunneth property that the map $CH_*(X)\otimes CH_*(Y)\to CH_*(X\times_k Y)$ is an isomorphism for all separated $k$-schemes $Y$, many classifying spaces $BG$ do satisfy this property. For example, Totaro proved\cite[Lemma 2.12]{burtbook} that, under certain conditions on the field $k$, the classifying space of any iterated wreath product of a finite abelian group satisfies the Chow Kunneth property. On the other hand, Totaro also gives examples\cite[Corollary 3.1]{burt} of finite groups $G$ such that $BG$ does not satisfy the Chow Kunneth property.

A somewhat stronger property than the Chow Kunneth property is the motivic Kunneth property, which gives information about all motivic homology groups of $X\times_k Y$, rather than just the Chow groups. 

In this paper, we turn from finite groups to the case of connected algebraic groups $G$. 

We show that if $G$ is any of the groups $GL(n),SL(n),O(q),Sp(2n), G_2$, or $SO(q)$ for an odd-dimensional quadratic form, then $BG$ satisfies the motivic Kunneth property over any field $k$ (of characteristic not equal to $2$ in the case of $O(q)$)[Propositions  \ref{Sp2n}, \ref{Oq} and \ref{G2}].

The even-dimensional special orthogonal group $SO(q)$ may or may not have the motivic Kunneth property, depending on the quadratic form $q$:
\\
\\\textbf{Proposition \ref{SOqeven}}: Let $k$ be an arbitrary field of characteristic not equal to $2$ and let $q$ be a nondegenerate quadratic form over $k$ of dimension $2n$. Then the classifying space of the special orthogonal group $SO(q)$ has the Chow Kunneth property if and only if  $\text{det}(q) = (-1)^n (\text{mod }  (k^\times)^2)$.
\\
\\In particular, the classifying space of the split form of $SO(2n)$ does have the Chow Kunneth property, even though some forms of $SO(2n)$ do not. One might conjecture that the classifying space of every split reductive group has the Chow Kunneth property.
\\
\\The technical tools for studying the Chow Kunneth property involve the triangulated category of motives $DM(k)$. Totaro defined the compactly supported motive of a quotient stack $M^c([X/G])\in DM(k)$ and proved that $[X/G]$ satisfies the motivic Kunneth property if and only if the motive $M^c([X/G])$ is a mixed Tate motive\cite{burt}. In the future, we will therefore substitute ``mixed-Tate property" for ``motivic Kunneth property".

In the first two sections of the paper we refine the construction of the $M^c([X/G])$, as background for the Chow Kunneth question. In section $1$, we construct the compactly supported motive $M^c(X)$ for any algebraic space $X$, and demonstrate that it satisfies expected properties, such as a localization exact triangle and representability of the Chow groups. In section $2$, we build upon the work of section $1$ to define $M^c([X/G])$ for any quotient stack $[X/G]$ where $X$ is an algebraic space and $G$ is an affine group scheme acting on it; this forms a slightly larger class of quotient stacks than considered by Totaro.

More importantly, we redefine $M^c([X/G])$ in a functorial way [Theorem \ref{mcfunctor}], since Totaro's construction in the triangulated category $DM(k)$ faces the issue that the cone in a triangulated category is only defined up to isomorphism, rather than canonical isomorphism. This functoriality enables a proof of a localization exact triangle for quotient stacks [Corollary \ref{localizationstacks}].

Functoriality of $M^c([X/G])$ has some concrete consequences. For example, it gives a simple proof that if $B(G\times H)$ is mixed Tate, then so are $BG$ and $BH$ [Lemma \ref{mixedtatefunctoriality}]. Moreover, the functoriality of $M^c([X/G])$ also aids in the proof of the following result:
\\
\\ \textbf{Proposition \ref{unipotent}}: Let $G$ be an affine group scheme and let $J\subset G$ be a normal split unipotent subgroup scheme of dimension $d$. Suppose $G/J$ acts on a scheme $Z$. Then $M^c([Z/G])\cong M^c([Z/(G/J)])(d)[2d]$. In particular, $M^c(BG)\cong M^c(B(G/J))(d)[2d]$.
\chapter{Defining motives of algebraic spaces}
Let $k$ be a field and let $R$ be a commutative ring. We work throughout the paper with the triangulated category $DM(k;R)$ of motives over $k$ with coefficients in $R$. In order to make use of the main properties of motives, we will assume throughout that the \textit{exponential characteristic} of $k$, defined to be $1$ if $\text{char} k = 0$ and $p=\text{char} k$ otherwise, is invertible in $R$. In this section, we work within the smaller triangulated subcategory $DM^{eff}_{-}(k;R)\subset DM(k;R)$ of ``effective, bounded above" motives, which is the setting first considered by Voevodsky.\cite{voe}

Given any separated scheme $X$ of finite type over $k$, Voevodsky defined the associated motives $M(X)$ as well as the compactly supported motive $M^c(X)$ in $DM^{eff}_{-}(k;R)$. In this section we extend Voevodsky's construction of $M^c(X)$ to the case of an algebraic space $X$ of finite type over $k$, and show that these motives satisfy the expected properties, such as functoriality, a localization exact triangle, as well as isomorphisms $Hom_{DM}(R(i)[2i],M^c(X))\cong CH_i(X;R)$. 

%Voevodsky constructed the triangulated category $DM^{eff}_{-}(k;R)$ of effective, bounded-below motives that provides a setting for motivic cohomology. Cisinski and Deglise constructed the larger category $DM(k;R)$, which contains $DM^{eff}_{-}(k;R)$ as a subcategory and in which the Tate motive is invertible. In addition, Voevodsky defined an assignment of a motive $M(X)\in DM^{eff}_{-}(k;R)$, as well as a ``compactly supported motive" $M^c(X)\in DM^{eff}_{-}(k;R)$, to any scheme $X$ of finite type over $k$. In this section, we extend Voevodsky's construction of compactly supported motives to the category of algebraic spaces. This is a straightforward extension of Voevodsky's definitions, using the fact that cycles in algebraic spaces behave in much the same way as they do in schemes. To show that our definition is natural, we prove a localization triangle and demonstrate that the Chow groups of algebraic spaces are representable in $DM(k;R)$. Along the way, we also define bivariant cohomology of algebraic spaces.

\section{Background}

We begin by reviewing the construction of $DM^{eff}_{-}(k;R)$ and its main properties, while at the same time generalizing relevant definitions to the context of algebraic spaces, which will be needed later. The main reference is \cite{voe}, as well as \cite{mvw}.

Throughout, all algebraic spaces will be assumed to be over a field $k$.
\begin{definition}
	\label{equi}
	A map $f:X\to S$ of algebraic spaces is equidimensional of dimension $r$ if it is locally of finite type, every irreducible component of $X$ is dominant over an irreducible component of $S$, and all the fibers are equidimensional of dimension $r$.\end{definition}

\begin{definition}
	\label{cor}
	Let $X$ be an algebraic space of finite type over $k$. An elementary correspondence from $S$ to $X$ is an integral subspace of $S\times X$. A correspondence from $S$ to $X$ is a cycle on $S\times X$.
	
	Define $z_{equi}(X,r)(S)$ to be the abelian group generated by elementary correspondences from $S$ to $X$ which are equidimensional of dimension $r$ over $S$. Define $\mathbb{Z}_{tr}(X)(S)=Cor(S,X)$ to be the subgroup of $z_{equi}(X,0)(S)$ generated by correspondences in which are finite over $S$. \end{definition}

Given a morphism $f:S\to X$, the cycle $[\Gamma_f]$ is a finite correspondence from $S$ to $X$, which we also denote by $f$.

We can compose correspondences between smooth $k$-schemes as follows.

Let $X,Y$ and $Z$ be smooth $k$-schemes, $\alpha$ a finite correspondence from $X$ to $Y$ and $\beta$ a finite correspondence from $Y$ to $Z$.

We use the notation $p^{XYZ}_{XY},p^{XYZ}_{XZ},p^{XYZ}_{YZ}$ for the projections $X\times Y\times Z$ to $X\times Y$, $X\times Z$, $Y\times Z$, respectively.

Then $\beta \circ \alpha$ os defined as $p^{XYZ}_{XZ*}(p_{YZ}^{XYZ*}\beta \cdot p^{XYZ*}_{XY}( \alpha))$. Voevodsky shows this is well-defined, composition satisfies transitivity and the diagonal $[\Delta]\subset X\times X$ acts as identity on $X$. This gives the set of finite correspondences of smooth $k$-schemes a category structure, denoted $Cor_k$, whose objects are smooth $k$-schemes and where $Mor(X,Y)=Cor(X,Y)$.

A \textit{presheaf with transfers} (with coefficients in $R$) is a functor $F:Cor_k\to R-Mod$. A  \textit{sheaf with transfers} (with respect to a given Grothendieck topology on smooth $k$-schemes) is a presheaf with transfers which is additionally a sheaf over the site of smooth $k$-schemes with the given topology.

The construction of $DM^{eff}_-(k;R)$ makes use of the Nisnevich Grothendieck topology. A family of morphisms of schemes $\{U_i\to X\}$ is a \textit{Nisnevich cover} if it is etale and for every point $x$ of $X$ there exists $i$ and a point $u$ of $U_i$ lying over $x$ such that the map $k(u)\to k(x)$ induced by $U_i\to X$ is an isomorphism. In particular, a Nisnevich cover is also an etale cover.

For any scheme $X$ of finite type over $k$, the presheaves $\mathbb{Z}_{tr}(X)$ and $z_{equi}(X,r)$ in fact have the structure of integral presheaves with transfer. In fact, they are sheaves with respect to the etale topology, hence in particular the Nisnevich topology.

If a coefficient ring $R$ is understood, and $X$ is a smooth scheme over $k$, we define the sheaves of $R$-modules with transfers $L(X)=R_{tr}(X):=\mathbb{Z}_{tr}(X)\otimes R$ and $L^c(X)=z_{equi}(X,0)\otimes R$.

Let $Sh_{Nis}(SmCor(k);R)$ denote the category of sheaves of $R$-modules with transfer on the site $Sm/k$ equipped with the Nisnevich topology. Let $\mathcal{C}:=Ch(Sh_{Nis}(SmCor(k);R))$ be the corresponding category of chain complexes. Let $D^- = D^-(Sh_{Nis}(SmCor(k);R))$ be its bounded-above derived category. A presheaf with transfers $F$ is \textit{homotopy invariant} if for every smooth $k$-scheme $X$, the projection $X\times \mathbb{A}^1\to X$ induces an isomorphism $F(X)\to F(X\times \mathbb{A}^1)$. Then $DM^{eff}_{-}(k;R)$ is the full triangulated subcategory of $D^{-}$ consisting of complexes with homotopy-invariant cohomology sheaves.

We may alternatively view $DM^{eff}_{-}(k;R)$ as a localization, rather than a subcategory, of $D^{-}$, as follows. Define the complex $C_*(F)$ for any presheaf with transfers $F$ by letting $C_i(F)(X)=F(X\times \Delta^i)$, where $\Delta^i$ is the algebraic simplex, and with differential given by the boundary map, the alternating sum of restriction to the faces. Then for any Nisnevich sheaf with transfers $F$, we have $C_*(F)\in DM^{eff}_{-}(k;R)$. The exact functor $C_*$ can be extended to an exact functor $C_*:D^{-}\to DM^{eff}_{-}(k;R)$ which is left-adjoint to the inclusion, and which identifies $DM^{eff}_{-}(k;R)$ as the localization of $D^{-}$ with respect to the localizing category generated by complexes of the form $L(X\times \mathbb{A}^1)\to L(X)$ induced by projection for any smooth $k$-scheme $X$.

For any scheme $X$ of finite-type over $k$, define $M(X)=C_*L(X)$ and $M^c(X)=C_*L^c(X)\in DM^{eff}_{-}(k;R)$. We will be primarily concerned with the latter, the compactly supported motive of $X$. If $X$ is proper over $k$, then $L^c(X)=L(X)$, so $M^c(X)=M(X)$. The motive $M(Spec(k))$ is denoted simply $R$. 

$DM^{eff}_{-}(k;R)$ is equipped with the structure of a tensor triangulated category, with $M(X)\otimes M(Y)=M(X\times_k Y)$.

An important motive is the Tate motive $R(1)$, with the property that the inclusion of a point $R=M(Spec(k))\to M(\mathbb{P}^1)$ is a direct summand whose complement is $R(1)[2]$. For any $r\ge 0$, the motive $R(r)$ denotes the $r$-fold tensor product $R(1)^{\otimes r}$ and, for any motive $M\in DM^{eff}_{-}(k)$, the motive $M(r)$ denotes the tensor product $M\otimes R(r)$.

While the functor $M(-)$ is covariantly functorial over arbitrary morphisms of schemes, $M^c(-)$ is only covariantly functorial for proper maps, and contravariantly functorial for flat maps, with dimension shift. More precisely, given a flat map $Y\to X$ of relative dimension $r$, there is a canonical pullback map $M^c(X)(r)[2r]\to M^c(Y)$.

Compactly supported motives satisfy homotopy invariance: given a vector bundle $E\to B$ of rank $r$, the flat pullback $M^c(B)(r)[2r]\to M^c(E)$ is an isomorphism in $DM^{eff}_{-}(k;R)$.

Given an open embedding of $k$-schemes $U\to X$ with closed complement $Z$, the natural sequence
\[
M^c(Z)\to M^c(X)\to M^c(U)
\]
extends to a distinguished triangle, the localization triangle. This last result was proved by Voevodsky in the case that $k$ has resolution of singularities, and extended to a general field by Kelly\cite{kelly}, under the assumption that the exponential characteristic of $k$ is invertible in $R$.

Finally, Chow groups are representable by compactly supported motives in $DM^{eff}_{-}(k;R)$, in the sense that we have natural isomorphisms

\[
Hom_{DM}(R(i)[2i], M^c(X))\cong CH_i(X)\otimes R
\]
\section{Defining motives and bivariant cohomology of algebraic spaces}
Let $X$ be an algebraic space of finite type over $k$. To generalize Voevodsky's definition of $M^c(-)$ for schemes, we work toward showing that the presheaves $z_{equi}(X,r)(S)$ defined in Definition \ref{cor} are in fact sheaves with transers.

We begin with some results about equidimensional cycles. The property of being equidimensional is not in general stable under pullback. For example, consider the morphism $f:Y\to X$ where $X\subset \mathbb{A}^2$ is given by $xy=0$, $Y\subset \mathbb{A}^3$ is given by $xz=0$, and the map $f:Y\to X$ is given by $(x,y,z)\mapsto (x,yz)$. This map is equidimensional of relative dimension $1$, but if we let $i:X_1\subset X$ be the closed subset $\{y=0\}\subset \mathbb{A}^2$, the pullback $f_1:Y_1\to X_1$ of $f$ by $i$ is not equidimensional. Indeed, one of the irreducible components of $Y_1$ is the line $\{x=0,y=0\}\subset \mathbb{A}^3$, which fails to be dominant over $X_1$. Therefore we need to distinguish between equidimensional maps and universally equidimensional maps.
\begin{lemma}
	\label{equilemma}
	a) Equidimensional maps of dimension $r$ are local on the domain in the etale topology. 
	\\b) Universally equidimensional maps of dimension $r$ are stable and local on the domain in the etale topology. 
	\\c) Let $f:X\to S$ be an equidimensional map of dimension $r$. If $S$ is a smooth $k$-scheme, then $f$ is universally open and universally equidimensional of dimension $r$.
\end{lemma}
\begin{proof}
	a) Let $f:X\to S$ be a morphism of algebraic spaces and let $u:X'\to X$ be an etale cover.
	
	Suppose $f$ is equidimensional of dimension $r$. Since etale maps are open, $u$ is open, hence each irreducible component of $X'$ dominates an irreducible component of $X$, hence also an irreducible component of $S$. Again by the fact that etale maps are open, for any $s\in S$, each component $W'$ of the fiber $X'_s$ dominates a component $W$ of the fiber $X_s$. Thus $dim W' = dim W = r$, so $f\circ u$ is equidimensional of dimension $r$.
	
	Conversely, suppose $f\circ u:X'\to S$ is equidimensional of dimension $r$. For each component $V$ of $X$, there exists by surjectivity of $u$ an irreducible component $V'$ of $X'$ mapping to $V$ in $X$. Since $f\circ u$ is equidimensional, we know that its restriction to $V'$ maps dominantly onto an irreducible component of $S$. But $V'\to S$ factors through $V$, hence $V$ also maps dominantly onto $S$. Similarly, for any $s\in S$, and any irreducible component $W$ of the fiber $X_s$, there is an irreducible component $W'$ of $X'_s$ mapping onto $W$. So, since the map $W'\to W$ has $0$-dimensional fibers, $dim W = dim W' = r$, demonstrating that $f$ is equidimensional of dimension $r$.
	\\
	\\b) Let $f:X\to S$ be a morphism, $u:S'\to S$ an etale cover, and $f':X'\to S'$ the base change of $f$ by $u$. Suppose $f'$ is equidimensional of dimension $r$. Then $u \circ f'$ is also equidimensional of dimension $r$. But $u\circ f'$ equals the map $X'\to X\to S$, and $X'\to X$ is an etale cover, so by part a), $f:X\to S$ is equidimensional of dimension $r$.
	\\
	\\c) When $X$ and $S$ are schemes and $f$ is of finite type, this is proved in \cite{voesus}, though this generalizes to $f$ locally of finite type since universally open and universally equidimensional maps are local in the domain. The claim then follows for $X$ and $S$ algebraic spaces since smoothness is stable, and universally open and universally equidimensional maps of dimension $r$ are stable and local in the domain.
\end{proof}
\noindent
Let $X$ and $Y$ be smooth schemes over $k$ and let $Z$ be an algebraic space of finite type over $k$. Let $\alpha$ be a finite correspondence from $X$ to $Y$ and let $\beta$ be a correspondence from $Y$ to $Z$ which is equidimensional of dimension $r$. Let $p^{XYZ}_{XY},p^{XYZ}_{XZ},p^{XYZ}_{YZ}$ denote the projections from $X\times Y\times Z$ to $X\times Y$, $X\times Z$, $Y\times Z$, respectively. We define 
\begin{equation*}
	{\beta\circ \alpha = p^{XYZ}_{XZ*}(p_{YZ}^{XYZ*}\beta \cdot_{p^{XYZ}_{XY}} \alpha)}
\end{equation*}.
Note: the above formula uses Fulton's refined intersection product\cite{fulton}, which can be evaluated along any morphism from a $k$-scheme to a smooth $k$-scheme. if $Z$ is smooth, then this definition coincides with that found in \cite{fulton}, where $\alpha\circ \beta$ is defined to be  $p^{XYZ}_{XZ*}(p_{YZ}^{XYZ*}\beta \cdot p^{XYZ*}_{XY}( \alpha))$
\\
\\In order for the above definition to give us a well-defined cycle and give the desired transfer maps, we need the following lemma:
\begin{lemma}
	\label{composition}
	In the above setting, 
	\\a) the cycle $p^{XYZ*}_{YZ}\beta$ intersects properly with $\alpha$ along $p^{XYZ}_{XY}$.
	\\b) The map $p^{XYZ}_{XZ}$ is proper when restricted to $|\alpha|\times Z \cap X\times |\beta|$.
	\\c) $\beta\circ \alpha$ is a correspondence from $X$ to $Z$ which is equidimensional of dimension $r$; if $r=0$ and $\beta$ is a finite correspondence, then $\beta\circ \alpha$ is also a finite correspondence.\end{lemma}
\begin{proof}We can assume that $\alpha$ and $\beta$ are elementary correspondences, so let $\alpha=[V]$ and $\beta=[W]$, for integral schemes $V\subset X\times Y, W\subset Y\times Z$, where $V$ is finite and surjective over $X$ and $W$ is equidimensional of dimension $r$ over $Y$. We can also assume $X$ and $Y$ are irreducible.
	
	a) We note that the intersection $V\times Z \cap X\times W \subset X\times Y\times Z$ equals the fiber product of $V\to Y$ with $W\to Y$, so denote it by $V\times_Y W$. Since $W$ is equidimensional of dimension $r$ over $Y$, then by Lemma \ref{equilemma}c), $V\times_Y W$ is equidimensional of dimension $r$ over $V$, so in particular all its components have dimension $dimV+r$. Since $dim W=dimY+r$, it follows that the codimension of each component of $V\times_Y W$ in $X\times W$ equals the codimension of $V$ in $X\times Y$, as desired.
	
	b) Note that $V\times_Y W$ is a closed subscheme of $V\times Z$, which is finite, hence proper, over $X\times Z$ since $V$ is finite over $X$.
	
	c) As noted above, $V\times_Y W$ is equidimensional of dimension $r$ over $V$, hence over $X$. Hence $p^{XYZ}_{XZ}(V\times_Y W)$ is also equidimensional of dimension $r$ over $X$, since the map $p^{XYZ}_{XZ}(V\times_Y W)\to X$ factors through $V\times_Y W\to X$ via the finite surjective map $V\times_Y W\to p^{XYZ}_{XZ}(V\times_Y W)$.
	
	Similarly, if $W$ is finite surjective over $Y$, then $V\times_Y W$ is finite surjective over $V$, hence over $X$, and all components of $p^{XYZ}_{XZ}(V\times_Y W)$ are also finite surjective over $X$.\end{proof}
\begin{lemma}
	\label{lci}
	Let $Y$ and $Y'$ be smooth algebraic spaces. Consider a fiber square
	\[
	\xymatrix{
		X' \ar[r]^{g'} \ar[d]^{f'} & Y' \ar[d]^{f}
		\\
		X \ar[r]^{g} & Y
	}
	\]
	Let $\alpha\in Z_*(Y'),\beta\in Z_*(X)$. Suppose that the map $X'\to X\times Y'$ is a regular closed embedding, and that $f':X'\to X$ as well as $f'\times id_{Y'}:X'\times Y'\to X\times Y'$ are l.c.i. morphisms. Suppose $f$ is proper along $|\alpha|$, the intersections $\beta \cdot_g f_*(\alpha)$, $f'^!\beta$ and $f'^!\beta \cdot_{g'} \alpha$ are proper intersections, and that $f'$ is proper along $|f'^!\beta \cdot_{g'} \alpha|$. Then $\beta \cdot_g f_*(\alpha) = f'_*(f'^!\beta \cdot_{g'} \alpha)$ in $Z_*(X)$.
\end{lemma}
\begin{proof}Assume $\alpha=[V], \beta=[W]$, where $V\subset Y', W\subset X$ are integral. Then consider the following fiber diagram:
	\[
	\xymatrix{
		f'^{-1}(W)\cap g'^{-1}(V) \ar[r] \ar[d]^{f'} & W\times V \ar[d]^{id_W\times f}
		\\
		W\cap g^{-1}(f(V)) \ar[r] \ar[d] & W\times f(V) \ar[d]
		\\
		X \ar[r]^{\gamma_g} & X\times Y
	}
	\]
	By definition, $[W]\cdot_g f_*[V]=\gamma_g^{!} ([W]\times f_*[V])$. By \cite[Theorem 6.2a)]{fulton}, we have $\gamma_g^{!} ([W]\times f_*[V])=f'_* \gamma_g^{!}([W]\times [V])$. On the other hand, we have a fiber diagram:
	\[
	\xymatrix{
		f'^{-1}(W)\cap g'^{-1}(V) \ar[r] \ar[d] & W\times V \ar[d]
		\\
		X' \ar[r]^{(f', g')} \ar[d] & X\times Y' \ar[d]
		\\
		X \ar[r]^{\gamma_g} & X\times Y
	}
	\]
	By \cite[Theorem 6.2c)]{fulton}, $\gamma_g^{!}([W]\times [V])=(f',g')^{!}([W]\times [V])$. Now, $(f',g')$ factors as $(f'\times id_{Y'})\circ \gamma_{g'}$. So $(f', g')^!([W]\times [V])=\gamma_{g'}^!((f'\times id_{Y'})^!([W]\times [V]))=\gamma_{g'}^!(f'^![W]\times[V])=f'^![W]\cdot_{g'} [V]$. Thus $[W]\cdot_g f_*[V]=f'_*(f'^![W]\cdot_{g'})$ as desired.
\end{proof}
\begin{lemma}
	\label{functoriality}
	Let $X,Y$ and $Z$ be smooth schemes and let $W$ be an algebraic space. Let $\alpha$ be a finite correspondence from $X$ to $Y$, $\beta$ a finite correspondence from $Y$ to $Z$, and $\gamma$ a quasifinite correspondence from $Z$ to $W$. Then 
	\\a) $\gamma \circ (\beta \circ \alpha) = (\gamma \circ \beta) \circ \alpha$.
	\\b) If $f:X\to Y$ is an arbitrary morphism, then $\beta \circ f$ equals the pullback $f^!\beta$, where $f^!:CH_*(|\beta|)\to CH_*(|(f\times id_Z)^{-1}(\beta)|)$ is the refined Gysin homomorphism corresponding to the following fiber square:
	\[
	\xymatrix{
		X\times Z \ar[r]^{f\times id_Z} \ar[d] & Y\times Z \ar[d]
		\\
		X \ar[r]^{f} & Y
	}
	\]
	\\c) If $f:X\to Y$ is a flat morphism, then $\beta \circ \Gamma_f$ equals the flat pullback $(f\times id_Z)^*(\beta)$.
\end{lemma}
\begin{proof}
	a) The proof is very similar to the proof found in \cite{fulton}, though it is modified to the case in which $W$ may not be smooth.
	We have
	\\
	$\gamma\circ (\beta\circ \alpha) = p^{XZW}_{XW*}(p^{XZW*}_{ZW}\gamma \cdot_{p^{XZW}_{ZW}} (\beta\circ \alpha))$
	\\
	$=p^{XZW}_{XW*}(p^{XZW*}_{ZW}\gamma \cdot_{p^{XZW}_{XZ}} p^{XYZ}_{XZ*}(p^{XYZ*}_{YZ}\beta \cdot_{p^{XYZ}_{XY}} \alpha))$
	\\
	$=p^{XZW}_{XW*}(p_{XZW*}(p^*_{XZ}\gamma \cdot_{p_{XYZ}} (p^{XYZ*}_{YZ}\beta \cdot_{p^{XYZ}_{XY}} \alpha)))$
	\\
	$=p_{XW*}((p^*_{ZW}\gamma \cdot_{p_{XYZ}} p^{XYZ*}_{YZ}\beta)\cdot_{p_{XY}}\alpha)$
	\\
	$=p_{XW*}((p^*_{ZW}\gamma \cdot_{p_{YZ}}\beta)\cdot_{p_{XY}}\alpha)$
	
	Where these equations follow by: i) and ii) The definition of composition; iii) Lemma \ref{lci}, which tells us in particular that, given cycles $\delta\in Z_*(X\times Z\times W), \epsilon\in Z_*(X\times Y\times Z)$, we have $\delta\cdot_{p^{XZW}_{XZ}} p^{XYZ}_{XZ*}\epsilon = p^*_{XZW}\delta \cdot_{p_{XYZ}} \epsilon$; iv) by composition of pushforwards and associativity of intersections; and v) by composition of intersection with pullback.
	
	On the other hand, we have:
	\\$(\gamma \circ \beta)\circ \alpha=p^{XYW}_{XW*}(p^{XYW*}_{YW}(\gamma \circ \beta)\cdot_{p^{XYW}_{XY}}\alpha)$
	\\$=p^{XYW}_{XW*}(p^{XYW*}_{YW}(p^{YZW}_{YW*}(p^{YZW*}_{ZW}\gamma \cdot_{p^{YZW}_{YZ}}\beta))\cdot_{p^{XYW}_{XY}}\alpha)$
	\\$=p^{XYW}_{XW*}(p_{XYW*}(p^*_{YZW}(p^{YZW*}_{ZW}\gamma\cdot_{p^{YZW}_{YZ}}\beta))\cdot_{p^{XYW}_{XY}}\alpha)$
	\\$=p^{XYW}_{XW*}(p_{XYW*}(p_{ZW}^*\gamma\cdot_{p_{YZ}}\beta)\cdot_{p_{XY}^{XYW}}\alpha)$
	\\
	$=p^{XYW}_{XW*}(p_{XYW*}((p^*_{ZW}\gamma\cdot_{p_{YZ}}\beta)\cdot_{p^*_{XY}}\alpha))$
	\\
	$=p_{XW*}((p^*_{ZW}\gamma \cdot_{p_{YZ}}\beta)\cdot_{p_{XY}}\alpha)$
	
	Where these equations follow by: i) and ii) the definition of composition; iii) Since $p^{XYW*}_{YW}\circ p^{YZW}_{YW*}=p_{XYW*}\circ p^*_{YZW}$; iv) composition of pullback with intersection; v) projection formula; vi) composition of pushforwards.
	\\
	\\b) Let $\gamma_f:X\to X\times Y$ denote the graph morphism. Then we have 
	\\$\beta \circ [\Gamma_f] = p^{XYZ}_{XZ*}(p_{YZ}^{XYZ*}\beta \cdot_{p^{XYZ}_{XY}} [\Gamma_f])$
	\\$=p^{XYZ}_{XZ*}(p_{YZ}^{XYZ*}\beta \cdot_{p^{XYZ}_{XY}} \gamma_{f*}[X])$
	\\$=p^{XYZ}_{XZ*}(\gamma_f\times id_Z)_*((\gamma_f\times id_Z)^!p_{YZ}^{XYZ*}\beta \cdot_{p^{XZ}_X} [X])$
	\\$=(\gamma_f\times id_Z)^!p_{YZ}^{XYZ*}\beta$
	\\$=f^!\beta$
	
	Where these equations follow by: i) the definition of composition; ii) $[\Gamma_f]=\gamma_{f*}[X]$; iii) Lemma \ref{lci}; iv) by composition of pushforwards, and the fact that intersecting with $[X]$ leaves a cycle unchanged; v) by composition of gysin homomorphisms, and the fiber diagram:
	\[
	\xymatrix{
		X\times Z \ar[r]^{\gamma_f\times id_z} \ar[d] & X\times Y\times Z \ar[r]^{p^{XYZ}_{YZ}} \ar[d] & Y\times Z \ar[d]
		\\
		X \ar[r]^{\gamma_f} & X\times Y \ar[r]^{p^{XY}_Y} & Y
	}
	\]
	\\
	c) Follows from b).
\end{proof}
\noindent It follows from Lemma \ref{functoriality} that $z_{equi}(X,r)$ and $\mathbb{Z}_{tr}(X)$ are presheaves with transfers. We now prove that they are in fact sheaves with transfers in the etale topology, hence in particular in the Nisnevich topology.
\begin{proposition}
	\label{sheaves}
	Let $X$ be an algebraic space of finite type over $k$. Then $z_{equi}(X,r)$ and $\mathbb{Z}_{tr}(X)$ are sheaves with transfers in the etale topology.\end{proposition}
\begin{proof}Let $S$ be a smooth scheme and let $\pi: S'\to S$ be an etale cover. Let $\pi_X:S'\times X\to S\times X$ denote $\pi\times id_X$.
	
	Identity: Let $\alpha,\beta\in z_{equi}(X,r)(S)$ and suppose $\alpha \circ \pi = \beta \circ \pi$. Thus, by Lemma \ref{functoriality}c), we have $\pi_X^*\alpha = \pi_X^* \beta$. In particular, $\pi_X^{-1}(|\alpha|)=|\pi_X^*\alpha|=|\pi_X^*\beta|=\pi_X^{-1}(|\beta|)$. So, since $h_{S\times X}$ is a sheaf in the etale topology, we have $|\alpha|=|\beta|$ as closed subschemes of $S\times X$. 
	
	So let $V_i$ be the irreducible components of $|\alpha|=|\beta|$. Let $W_{ij}\subset S'\times X$ be the irreducible components of $\pi_X^{-1}(V_i)$. Note that the $W_{ij}$, ranging over all $i$ and $j$, are distinct. Thus, in order for $\pi_X^*\alpha$ to equal $\pi_X^*\beta$, $\alpha$ and $\beta$ must have equal coefficients for all the $[V_i]$, so $\alpha=\beta$.
	
	Gluability: Let $\alpha'\in z_{equi}(X,r)(S')$. Let $p_1,p_2: S'\times_S S'\to S'$ denote the two projection maps and let ${p_{1}}_X, {p_{2}}_X: S'\times_S S'\times X\to S'\times X$. Note that ${p_{1}}_X, {p_{2}}_X$ can also be viewed as the projections $(S'\times X) \times_{S\times X} (S'\times X)\to S'\times X$. Suppose $\alpha'$ satisfies $\alpha'\circ p_1 = \alpha'\circ p_2$. So by Lemma \ref{functoriality}c), we have ${p_1}_X^*\alpha' = {p_2}_X^*\alpha'$. In particular, ${p_1}_X^{-1}(|\alpha'|) = {p_2}_X^{-1}(|\alpha'|)$, so, since $h_{S\times X}$ is a sheaf in the etale topology, and closed immersions are stable in the etale topology, there is a closed reduced subscheme $V\subset S\times X$ such that $\pi_X^{-1}(V)=|\alpha'|$. 
	
	Then we want to show that there exists a cycle on $V$ whose pullback by $\pi_X$ equals $|\alpha'|$. We may reduce to the case of $V$ irreducible. Let $V'_{i}$ be the irreducible components of $\pi_X^{-1}(V)$. So the $V'_{i}$ are distinct and are the irreducible components of $|\alpha'|$. Let $V''_{ij}=V'_{i}\times_V V'_{j}\subset (S'\times X)\times_{S\times X} (S'\times X)$. Again, the irreducible components of $V''_{ij}$, ranging over $i,j$, are distinct. Let $\alpha'=\sum_{i}n_{i}[V'_{i}]$. Then ${p_1}_X^*\alpha' = \sum_{i,j} n_i[V''_{i,j}]$, whereas ${p_2}_X^*\alpha'=\sum_{i,j} n_j[V''_{i,j}]$. thus, $n_i=n_j=n$ for all $i,j$. So indeed $\alpha'=f^*(n[V])$.
	
	It remains to show that, given a cycle $\alpha\in Z_*(S\times X)$, $\alpha\in z_{equi}(X,r)(S)$ if and only if $\alpha':=\alpha\circ \pi\in z_{equi}(X,r)(S')$ and $\alpha\in L(X)(S)$ if and only if $\alpha'\in L(X)(S')$. This is because the class of equidimensional maps of dimension $r$ is stable in the etale topology, as is the class of finite maps.
\end{proof}
Tensoring with $R$, we obtain Nisnevich sheaves of $R$-modules with transfers, $L^c(X)=z_{equi}(X,0)\otimes R$ and $L(X)=R_{tr}(X):=\mathbb{Z}_{tr}(X)\otimes R$. We can now define the motives $M^c(X)$ and $M(X)$.
%\\
%\\The following is proved in \cite{voe}.
%\begin{lemma}
%\label{homotopyinvariance}
%For any Nisnevich sheaf with transfers $F$, the complex $C_*(F)$ of Nisnevich sheaves with transfers has homotopy-invariant cohomology sheaves, hence gives an element of $DM^{eff}_{-}(k;R)$.
%\end{lemma}
\begin{definition}
	\label{motive}
	Let $X$ be an algebraic space of finite type over $k$. Then we define $M^c(X):=C_*(L^c(X))$ and $M(X)=C_*(L(X))$\end{definition}
\noindent We also define the bivariant cohomology groups $A_{r,i}(Y,X)$, though for simplicity we restrict ourselves to the case where $Y$ is a smooth $k$-scheme.
\begin{definition}
	\label{bivariant}
	Let $X$ be an algebraic space of finite type over $k$ and let $Y$ be a smooth $k$-scheme. Let $r\ge 0$ and $i$ be integers. Then we define $A_{r,i}(Y,X)$ by the hypercohomology $\textbf{H}_{Nis}(Y,C_*(z_{equi}(X,r)))$. Let $A_{r,i}(X):=A_{r,i}(Spec k, X)$.\end{definition}
In the case of schemes, the groups $A_{r,i}(X)$ are the Borel-Moore motivic homology groups:
\[
A_{r,i}(X)\cong Hom_{DM(k;R)}(R(r)[2r+i],M^c(X))=H^M_{2r+i}(X,R(r))
\]
In particular, $A_{r,0}(X)\cong CH_r(X)$.

Proposition \ref{bivariance} provides another interpretation of the bivariant cohomology groups $A_{r,i}(Y,X)$.
\section{Localization}
We'll make use of the following theorem from 
%\cite{voe} and from 
\cite{gruray}:
%\begin{theorem}[Voevodsky]
%\label{birational}
%Let $F$ a presheaf with transfers on $Sm/k$ such that for any smooth scheme $S$ over $k$ and a section $\phi\in F(S)$, there is a proper birational morphism $p:S'\to S$ with $F(p)(\phi)=0$. Then the complex $C_*(F)$ is quasiisomorphic to $0$.\end{theorem}
\begin{theorem}(platification\cite{gruray})
	\label{platification}
	Let $S$ be a quasicompact and quasiseparated algebraic space, $U\subset S$ an open subspace, $f:X\to S$ an algebraic space of finite presentation over $S$, $\mathcal{M}$ a finite $\mathcal{O}_X$-module such that $\mathcal{M}|_f^{-1}(U)$ is flat over $U$. Then there is a blowup $p:S'\to S$ along a closed subscheme that lies outside of $U$ such that the strict transform of $\mathcal{M}$ along $p$ is flat over $S'$.\end{theorem}
\begin{lemma} (generic flatness for algebraic spaces)
	\label{genericflatness}
	Let $X$ be an algebraic space of finite type over $k$ and $S$ a reduced scheme of finite type over $k$. Let $f:X\to S$ be a morphism. Then there exists a dense open subset $U\subset S$ such that $X_U$ is flat over $U$. \end{lemma}
\begin{proof}The case of $X$ a scheme is proved (in greater generality) in \cite[Tag 0529]{stacks}. If $X$ is an algebraic space, then since $X$ is quasicompact, there is a cover $X'\to X$ of $X$ by an etale map with $X'$ quasicompact. So in particular $X'$ is of finite type over $k$. Applying generic flatness to $X'\to S$, let $U$ be a dense open set of $S$ such that $X'_U$ is flat over $U$. Therefore, $X_U$ is flat over $U$, since $X'_U$ is an etale cover of $X_U$ and flatness is etale-local in the domain.
\end{proof}
\begin{lemma}
	\label{equilemmageneric}
	Let $X$ be an algebraic space, $S$ a scheme and let $f:X\to S$ a flat morphism of finite type. Then $f$ is equidimensional of dimension $r$ if and only if for all generic points $y: Spec(K)\to S$ of the irreducible components of $S$, the projection $X\times_S Spec(k)\to Spec(K)$ is equidimensional of dimension $r$.\end{lemma}
\begin{proof}
	In the case that $X$ is a scheme, this is proven in \cite{voesus}. It then follows for $X$ an algebraic space by the fact that equidimensional maps of dimension $r$ are etale-local in the domain.
\end{proof}
\begin{proposition}
	\label{prelocalization}
	Let $X$ be an algebraic space of finite type over $k$, where $k$ admits resolution of singularities. Let $Z$ be a closed subspace of $X$ and let $U:=X-Z$. Then for any $r\ge 0$ there is a natural left-exact sequence of complexes of Nisnevich sheaves with transfers
	\begin{equation*}
		0\to C_*(z_{equi}(Z,r)\otimes R)\to C_*(z_{equi}(X,r)\otimes R)\to C_*(z_{equi}(U,r)\otimes R) 
	\end{equation*}
	and $coker(C_*(z_{equi}(X,r)\otimes R)\to C_*(z_{equi}(U,r)\otimes R))$ is acyclic.
\end{proposition}
\begin{proof}
	Since $C_*$ is an exact functor from presheaves with transfers to complexes of presheaves with transfers, it suffices to prove that there is a natural left-exact sequence of presheaves
	\begin{equation*}
		0\to z_{equi}(Z,r)\otimes R\to z_{equi}(X,r)\otimes R\to z_{equi}(U,r)\otimes R
	\end{equation*}
	such that the presheaf cokernel $Q:= coker(z_{equi}(X,r)\to z_{equi}(U,r))$ has $C_*Q$ acyclic. 
	
	Indeed, the exactness of the above sequence follows directly from the definitions, so it remains to show that $C_*Q$ is acyclic. By \cite[Theorem 5.3.1]{kelly}, it suffices to show that, viewing $Q$ as a presheaf over $Sch/k$ rather than just $SmSch/k$, for any scheme $S\in Sch/k$ and section $\phi\in F(S)$, there is a proper birational morphism $p:S'\to S$ with $F(p)(\phi)=0\in F(S')$.
	
	Indeed, let $\phi\in Q(S)$ and let $\alpha\in z_{equi}(U,r)(S)$ be a representative of $\phi$. Let $W|_U\subset S\times U$ denote the support of $\alpha$ and let $W\subset S\times X$ be the closure of $W|_U$ in $S\times X$. Let $f:W\to S$ denote the projection. By generic flatness, there is a dense open subset $O\subset S$ such that $f^{-1}(O)$ is flat over $O$. Thus by platification, there exists a blowup map $p:S'\to S$ along a closed subscheme that lies outside of $O$ such that the strict transform of $W$, which we denote $W'$, is flat over $S'$. The map $p$ is proper and birational.
	
	Let $W'|_U$ denote $W\cap (S'\times U)$, which is also the same as the strict transform of $W|_U$ along $p$. We claim that, in fact, $W'|_U$ equals $S'\times_S W|_U$, the total transform of $W|_U$. For this we need to show that $E\times_S W|_U$ does not contain any irreducible component of $S'\times_S W|_U$, where $E\subset S'$ is the exceptional divisor. Indeed, this is true by dimension reasons: since $W|_U$ is equidimensional of dimension $r$ over $S$, Lemma \ref{equilemma}c) implies that all components of $S'\times_S W|_U$ are equidimensional of dimension $r$ over $S'$ (which is irreducible by our assumption that $S$ is irreducible), thus of dimension $dimS' + r$. But by the same token, all components of $E\times_S W|_U$ are equidimensional of dimension $r$ over $E$, hence of dimension $dimE + r$. Since $dimE=dimS'-1$, we conclude that $E\times_S W|_U$ contains no irreducible component of $S'\times_S W|_U$. So $W'\cap (S'\times U) = |\alpha\circ p|$. Let $\beta$ be a cycle on $W'$ such that the restriction of $\beta$ to $W'|_U$ equals $\alpha \circ  p$. Note $W'|_U$ is dense in $W'$ since $W|_U$ is dense in $W$. So $W'$ is generically equidimensional of dimension $r$ over $S'$, hence by flatness and Lemma \ref{equilemmageneric}, it is equidimensional of dimension $r$ over $S'$. Thus $\beta\in z_{equi}(X,r)(S')$. So $ p: S'\to S$ is the desired proper birational morphism such that $Q(q\circ p)(\phi)=0$.
\end{proof}
\begin{corollary}
	\label{localization}
	In the setting of Proposition \ref{prelocalization}, there is a canonical localization triangle of compactly supported motives:
	\begin{equation*}
		M^c(Z)\to M^c(X)\to M^c(U)\to M^c(Z)[1]
	\end{equation*}
\end{corollary}
\begin{corollary}
	\label{bivariantseq}
	Let $X,Z$ and $U$ be as in the setting of Proposition \ref{prelocalization}, let $Y$ be a smooth $k$-scheme and let $r\ge 0$. Then there is a canonical long exact sequence of bivariant cohomology groups:
	\begin{equation*}
		\cdots \to A_{r,i+1}(Y,U)\to A_{r,i}(Y,Z)\to A_{r,i}(Y,X)\to A_{r,i}(Y,U)\to A_{r,i-1}(Y,Z)\to\cdots
	\end{equation*} 
\end{corollary}
\noindent Totaro defined\cite{burt}, for any $j\in \mathbb{Z}$, the subcategory $D_j(k)\subset DM(k;R)$ as the smallest localizing subcategory containing $M^c(X)(a)$ for separated $k$-schemes $X$ and integers $a$ such that $\text{dim} X+a\le j$. An analogous filtration was considered by Voevodsky\cite{voe} for effective motives. Then this subcategory also captures the dimension of algebraic spaces, in the following sense:
\begin{corollary}
	\label{filtration}
	Let $X$ be an algebraic space of finite type over $k$ and let $a\in \mathbb{Z}$ such that $\text{dim} X + a\le j$. Then $M^c(X)(a)\in D_j(k)$.
\end{corollary}
\begin{proof}
	Let $U\subset X$ be an open subset with complement $Z\subset X$ such that $U$ is a scheme and $Z$ has lower dimension than $X$. Then the corollary follows by the localization triangle and induction.
\end{proof}
\section{Representability of bivariant cohomology and Chow groups}
\begin{proposition}
	\label{chowbivariant}
	Let $X$ be an algebraic space of finite type over $k$ and let $r\ge 0$. Then $A_{r,0}(X)=CH_r(X)\otimes R$.\end{proposition}
\begin{proof} The hypercohomology of $Spec k$ with coefficients in the complex of sheaves $C_*(z_{equi}(X,r))$ is simply the cohomology of the complex $C_*(z_{equi}(X,r))(Spec k)$. In particular, $A_{r,0}(X)$ is the cokernel of the map \[
	z_{equi}(X,r)(\textbf{A}^1)\otimes R\to z_{equi}(X,r)(Spec k)\otimes R\] given by $i_1^*-i_0^*$, where $i_0,i_1: Speck\to \textbf{A}^1$ correspond to the points $0$ and $1$, respectively, in $\textbf{A}^1$. By right-exactness of tensoring, we may reduce to the case $R=\mathbb{Z}$.
	
	But note that an elementary cycle $[V]\in Z_{r+1}(\textbf{A}^1\times X)$ is dominant over $\textbf{A}^1$ if and only if it is equidimensional of dimension $r$ over $\textbf{A}^1$. Indeed, the irreducible components of any fiber have dimension at least $r$; if one of the fibers had dimension $r+1$, then it would be all of $V$, so $V\to \textbf{A}^1$ would be constant, hence not dominant. Thus $A_{r,0}(X)= CH_r(X)$ as desired.\end{proof}
\noindent The following theorem is proved in \cite{voe}.
\begin{theorem}
	\label{representabilityschemes}
	Let $Y$ be a smooth $k$-scheme and $F$ a presheaf with transfers on $Sm/k$. Then for any $i\ge 0$, there are canonical isomorphisms \begin{equation*}Hom_{DM^{eff}_{-}(k;R)}(M(Y),C_*(F)[i])\cong \textbf{H}^i_{Nis}(Y,C_*(F))\end{equation*} and these isomorphisms are natural in $F$.\end{theorem}
\begin{corollary}
	\label{representabilitycor}
	Let $Y$ be a smooth $k$-scheme and $X$ an algebraic space of finite type over $k$. Then there are canonical isomorphisms $Hom_{DM}(M(Y)[i],M^c(X))\cong A_{0,i}(Y,X)$.\end{corollary}
\noindent We construct canonical morphisms \begin{equation*}\varphi_{r,i}:Hom_{DM^{eff}_{-}(k;R)}((M(Y)(r)[2r+i], M^c(X))\to A_{r,i}(Y,X)\end{equation*} 
Corollary \ref{representabilitycor} gives us canonical isomorphisms \begin{equation*}\varphi_{0,i}:Hom_{DM^{eff}_{-}(k;R)}((M(Y)[i], M^c(X))\to A_{0,i}(X)\end{equation*}
Now suppose $r>0$. Then $C_*(Y)(r)[2r]$ is canonically isomorphic to the kernel of the projector on $C_*(Y\times \textbf{P}^1)(r-1)[2r-2]$ induced by the composition $\textbf{P}^1\to Speck\to \textbf{P}^1$. Thus, we have a canonical map \begin{equation*}Hom_{DM^{eff}_{-}(k;R)}((M(Y)(r)[2r+i], M^c(X))\to\end{equation*} \begin{equation*}Hom_{DM^{eff}_{-}(k;R)}((M(Y\times \textbf{P}^1)(r-1)[2r-2+i], M^c(X))\end{equation*}
By induction, there is a canonical map \begin{equation*}\varphi_{r-1,i}:Hom_{DM^{eff}_{-}(k;R)}((M(Y\times \textbf{P}^1)(r-1)[2r-2+i], M^c(X))\to A_{r-1,i}(Y\times \textbf{P}^1,X)\end{equation*}
We also have a canonical map $A_{r-1,i}(Y\times \textbf{P}^1,X)\to A_{r,i}(Y,X\times \textbf{P}^1)$, given by the natural morphisms $z_{equi}(X,r-1)(Y\times \textbf{P}^1\times \Delta^n)\to z_{equi}(X\times \textbf{P}^1,r)(Y)$. Finally, there is a map $A_{r,i}(Y,X\times \textbf{P}^1)\to A_{r,i}(Y,X)$ by pushing forward along the projection $X\times \textbf{P}^1\to X$. We then define \begin{equation*}\varphi_{r,i}:Hom_{DM^{eff}_{-}(k;R)}((M(Y)(r)[2r+i], M^c(X))\to A_{r,i}(Y,X)\end{equation*}
to be the composition of the above maps.
\\
\\The following theorem is proved in \cite{voe}.
\begin{theorem}
	\label{bivarianceschemes}
	If $X$ is a scheme of finite type over $k$, the maps \begin{equation*}\varphi_{r,i}:Hom_{DM^{eff}_{-}(k;R)}((M(Y)(r)[2r+i], M^c(X))\to A_{r,i}(Y,X)\end{equation*}defined above are isomorphisms for all $i$ and all $r\ge 0$.\end{theorem}
\begin{proposition}
	\label{bivariance}
	Let $X$ be an algebraic space of finite type over $k$ and let $Y$ be a smooth $k$-scheme. Then the maps \begin{equation*}\varphi_{r,i}:Hom_{DM^{eff}_{-}(k;R)}((M(Y)(r)[2r+i], M^c(X))\to A_{r,i}(Y,X)\end{equation*}defined above are isomorphisms for all $i$ and all $r\ge 0$.\end{proposition}
\begin{proof} For convenience let $B_{r,i}(Y,X)$ denote $Hom_{DM^{eff}_{-}(k;R)}((M(Y)(r)[2r+i], M^c(X))$. 
	
	Let $U$ be a dense open subscheme of $X$, and $Z=X-U$. By localization, the triangle $M^c(Z)\to M^c(X)\to M^c(U)$ is a distinguished triangle, so it induces a long exact sequence \begin{equation*}\cdots \to B_{r,i+1}(Y,U)\to B_{r,i}(Y,Z)\to B_{r,i}(Y,X)\to B_{r,i}(Y,U)\to B_{r,i-1}(Y,Z)\to\cdots\end{equation*}
	
	By the localization sequence for bivariant cohomology, we also have a long exact sequence \begin{equation*}\cdots \to A_{r,i+1}(Y,U)\to A_{r,i}(Y,Z)\to A_{r,i}(Y,X)\to A_{r,i}(Y,U)\to A_{r,i-1}(Y,Z)\to\cdots\end{equation*}
	
	Since the morphisms $\varphi_{r,i}:B_{r,i}(Y,-)\to A_{r,i}(Y,-)$ are natural, they give a morphism of complexes between these two long exact sequences. We proceed by induction on the dimension of $X$ (where the base case of $0$-dimensional spaces holds since such spaces are necessarily schemes), so we can assume the result of the proposition holds for $Z$. By Theorem \ref{bivarianceschemes}, the result is also true for the scheme $U$. Thus by the five lemma, the result holds for $X$.\end{proof}
\noindent Applying Proposition \ref{bivariance} to the special case of $Y=Spec k$, $i=0$ gives us the following corollary.
\begin{corollary}
	\label{chowrep}
	Let $X$ be an algebraic space of finite type over $k$. Then there are canonical isomorphisms \begin{equation*}\varphi_{r}:Hom_{DM^{eff}_{-}(k;R)}((\textbf{Z}(r)[2r], M^c(X))\to CH_r(X)\otimes R\end{equation*}for all $r\ge 0$.
	See \cite{edidingraham} for a discussion of the Chow groups of algebraic spaces.
\end{corollary}
\chapter{A functorial construction of compactly supported motives of quotient stacks}
We recall Totaro's construction\cite{burt} of the compactly supported motive of a quotient stack, while slightly generalizing it. Let $\mathcal{Q}$ be the category of quotient stacks of the form $[Y/G]$, where $Y$ is any algebraic space of finite type over $k$ and $G$ is an affine group scheme acting on it over $k$. Let $X\in \mathcal{Q}$.

Define a \textit{resolution of $X$ by algebraic spaces} to be the data of a sequence of surjections $\cdots\xrightarrow{f_2} V_2\xrightarrow{f_1} V_1$ of vector bundles over $X$ and closed substacks $S_i\subset V_i$ such that $V_i-S_i$ is an algebraic space and $S_{i+1}\subset f_i^{-1}(S_i)$, and such that the codimension of $S_i$ in $V_i$ approaches infinity. For simplicity, refer to this data as $(V_i-S_i)$. Such a resolution induces a zig-zag of algebraic spaces
\[
\cdots \to V_3-S_3 \xleftarrow{i_2} V_3-f_2^{-1}(S_2)\xrightarrow{f_2} V_2-S_2 \xleftarrow{i_1} V_2-f_1^{-1}(S_1)\xrightarrow{f_1} V_1-S_1
\]
Where the $f_k$ are vector bundles and the $i_k$ are open embeddings. Since pulling back by a vector bundle induces an isomorphism on compactly supported motives, we can invert these pullbacks.
\begin{definition}
	\label{mcx}
	Define $M^c(X)$ as the homotopy limit of the sequence of motives
	\[
	\cdots \xrightarrow{(f_2^*)^{-1}} M^c(V_2-S_2) \xrightarrow{i_1^*} M^c(V_2-f_1^{-1}(S_1))\xrightarrow{(f_1^*)^{-1}} M^c(V_1-S_1)
	\]
	for any resolution of $X$ by algebraic spaces $(V_i-S_i)$.
\end{definition}
In fact, Totaro's definition\cite{burt} applied only to the case $X=[Y/G]$, where $Y$ is a quasi-projective $k$-scheme with a $G$-equivariant ample line bundle; these assumptions were necessary to guarantee $X$ has a resolution by schemes. However, since we can now make use of compactly supported motives of algebraic spaces, we may use a resolution by algebraic spaces rather than schemes, and therefore no longer need these assumptions. The key is that we may find by Totaro\cite{burtchow} a sequence of representations $V_i$ of $G$ and closed subsets $S_i\subset V_i$ of codimensions approaching infinity such that $G$ acts freely on $V_i-S_i$. But the quotient of any scheme by a free group action is an algebraic space, so $[Y/G]$ has a resolution by the algebraic spaces $((V_i-S_i)\times Y)/G$.

In the context considered in his paper, Totaro proved\cite{burt} that the motive $M^c(X)$ does not depend, up to isomorphism, on the choice of resolution $(V_i-S_i)$. The same proof works in the larger context considered here. However, homotopy limits in triangulated categories are defined only up to isomorphism, not canonical isomorphism, so $M^c$ is not clearly functorial. The goal of this section is to give a functorial construction of $M^c$.

We generalize the problem a little by considering homotopy limits of generalized zig-zags. Let $\mathcal{A}$ be the category whose objects are diagrams of algebraic spaces of finite type over $k$ of the following form:
\[
\cdots \to X_3 \xleftarrow{i_2} Y_2 \xrightarrow{j_2} E_2\xrightarrow{f_2} X_2 \xleftarrow{i_1} Y_1\xrightarrow{j_1} E_1\xrightarrow{f_1} X_1
\]
where $f_k$ are vector bundles, $j_k$ are closed embeddings and $i_k$ are open embeddings. We denote this object by $(X_i)$, omitting the data of the $Y_i$ and $E_i$ from the name.

For most objects of $\mathcal{A}$ that we look at, the closed embeddings $j_k$ will be the identity maps, though we will occasionally need them to be proper closed embeddings. As an example, given a quotient stack $X$ over $k$, a resolution of $X$ by algebraic spaces, as defined above, can be viewed as an object of $\mathcal{A}$.

Define a morphism in $\mathcal{A}$ to be a morphism of diagrams between the given zig-zags. 

In practice, however, we focus on two classes of morphisms in $\mathcal{A}$. A \text{proper morphism} in $\mathcal{A}$ from $(X_i)$ to $(X'_i)$ is a morphism of diagrams
\[
\xymatrix{
	\cdots \ar[r] & X_2 \ar[d] & Y_1 \ar[d] \ar[l] \ar[r] & E_1 \ar[r] \ar[d] & X_1 \ar[d] \\
	\cdots \ar[r] & X'_2 & Y'_1 \ar[l] \ar[r] & E'_1 \ar[r] & X'_1
}
\]
where the vertical maps are all proper, and the squares
\[
\xymatrix{
	E_i \ar[r] \ar[d] & X_i \ar[d]\\
	E'_i \ar[r] & X'_i
}, 
\xymatrix{
	X_{i+1} \ar[d] & Y_i \ar[l] \ar[d]\\
	X'_{i+1} & Y'_i \ar[l]
}
\]
are fiber squares.

On the other hand, a \textit{flat map} from $(X_i)$ to $(X'_i)$ is a morphism of diagrams as above, where the vertical maps are flat and the squares 
\[
\xymatrix{
	Y_i\ar[r] \ar[d] & E_i \ar[d]\\
	Y'_i \ar[r] & E'_i
}
\]
are fiber squares. The relative dimension of such a flat map is defined to be the relative dimension of the map $X_1\to X'_1$.

In summary, whether for proper maps or flat maps, we require all squares which mix flat and proper maps to be fiber squares.
\\
\\Let $\mathbb{N}^{op}$ denote the diagram category
\[
\cdots \to \cdot \to \cdot
\]
\begin{definition}
	\label{McA}
	We extend $M^c:AlgSp/k\to DM(k)$ to a functor $M^c:\mathcal{A}\to DM(k)^{\mathbb{N}^{op}}$ as follows: given an object $X=(X_i)\in \mathcal{A}$, define $M^c(X)$ to be the sequence
	\[ 
	\cdots \to M^c(X_2)(-d_1)[-2d_1]\to M^c(X_1)
	\]
	where $d_k$ is the rank of $E_k$ as a vector bundle over $X_k$, and the map 
	\[M^c(X_{k+1})(-d_k-\cdots - d_1)[-2d_k-\cdots - 2d_1]\to M^c(X_k)(-d_{k-1} - \cdots - d_1)[-2d_{k-1}-\cdots -2d_1]\]
	is the composition of open pullback induced by $i_k$, proper pushforward induced by $j_k$, and the inverse of the flat pullback induced by $f_k$.
	\\
	\\We then define the limit motive $LM^c(X)\in DM(k)$ to be the homotopy limit of the sequence $M^c(X)$, for any $X\in \mathcal{A}$.
\end{definition}
\noindent In this new terminology, the compactly supported motive $M^c(X)$ of a quotient stack $X\in \mathcal{Q}$ is $LM^c((V_i-S_i))$, for any resolution of $X$ by algebraic spaces $(V_i-S_i)$.

Just as on $AlgSp/k$, the construction $M^c$ is covariantly functorial on $\mathcal{A}$ with respect to proper maps, and contravariantly functorial (with dimension shift) with respect to flat maps.

However, as mentioned before, the homotopy limit $LM^c\colon \mathcal{A}\to DM(k)$ is not obviously functorial. Our goal is to make the construction $LM^c\colon \mathcal{A}\to DM(k)$ functorial (with respect to either flat maps or proper maps) and then to show that this functor descends to a functor on quotient stacks. In order to do this, we carry out the construction in a model category, rather than in the triangulated category $DM(k)$.

\section{Background}
We briefly recall the necessary background about the model category. Our main reference for this section is \cite{cisdeg}. Let $k$ be a field and let $SmCor(k)$ denote the category of smooth schemes, with finite correspondences as morphisms. Let $Sh_{Nis}(SmCor(k))$ denote the category of sheaves on the site $SmCor(k)$ equipped with the Nisnevich topology and let $\mathcal{C}:=Ch(Sh_{Nis}(SmCor(k)))$. To each scheme $X$ of finite type over $k$ (or, simply, a $k$-scheme) we may associate an object $\mathbb{Z}_{tr}(X)\in \mathcal{C}$, the Nisnevich sheaf with transfers represented by $X$. 

The category $\mathcal{C}$ can be equipped with a monoidal model category structure in which weak equivalences are $\mathbb{A}^1$-weak equivalences and cofibrations are maps in $I$-cof, where $I$ is the set of maps of the form $S^{n+1} \mathbb{Z}_{tr}(X)\to D^n \mathbb{Z}_{tr}(X)$, for $n\in \mathbb{N}$ and $X\in Sm/k$. In particular, for any $k$-scheme $X$, the object $\mathbb{Z}_{tr}(X)$ is cofibrant. The homotopy category $\text{Ho}\mathcal{C}$ is $DM^{eff}(k)$.
\\
\\We can now enlarge the model category $\mathcal{C}$ by adding spectra. Given any cofibrant object $S\in \mathcal{C}$, we define a symmetric $S$-spectrum to be a sequence of objects $(E_n)_{n\ge 0}$, where $E_n\in \mathcal{C}$, together with maps $\sigma_n: S\otimes E_n\to E_{n+1}$, such that the symmetric group $\Sigma_n$ acts on $E_n$ (via automorphisms in $\mathcal{C}$), and such that the composition
\begin{equation*}
	S^{\otimes m}\otimes E_n \to \cdots \to E_{m+n}
\end{equation*}
induced by the maps $\sigma_i$ is $\Sigma_m\times \Sigma_n$-equivariant.

We define a morphism between $S$-spectra $(E_n,\sigma_n)$ and $(F_n,\tau_n)$ to consist of a sequence of $\Sigma_n$-equivariant maps $f_n:E_n\to F_n$ which also respect the maps $\sigma_n, \tau_n$. We denote the category of symmetric $S$-spectra in $\mathcal{C}$ by $Sp_S(\mathcal{C})$. 

There is a monoidal stable model category structure on $Sp_S(\mathcal{C})$ where weak equivalences are stable $\mathbb{A}^1$-equivalence. There is a natural monoidal left-Quillen functor $\Sigma^\infty:\mathcal{C}\to Sp_S(\mathcal{C})$ given by $\Sigma^\infty(X)_n=X\otimes S^{\otimes n}$.

If $X\in \mathcal{C}$ is cofibrant, then $\Sigma^\infty(X)$ is cofibrant, so $\Sigma^\infty(X)\otimes -$ is left-Quillen, with adjoint denoted $\underline{Hom}_S(-,\Sigma^\infty(X))$. The essential property of $Sp_S^\Sigma(\mathcal{C})$ is that the functor $\Sigma^\infty(S)\otimes -:Sp_S^\Sigma(\mathcal{C})\to Sp_S^\Sigma(\mathcal{C})$ is a Quillen equivalence. We'll denote the functor $\Sigma^\infty(S)\otimes -$ by $\Sigma_S$ and its adjoint by $\Omega_S$. 
\begin{proposition}
	\label{spectra}
	Let $T,S\in \mathcal{C}$ be cofibrant objects. Suppose that they become isomorphic in the homotopy category $Ho(\mathcal{C})$. Then the homotopy categories $Ho(Sp_S^\Sigma(\mathcal{C}))$ and $Ho(Sp_T^\Sigma(\mathcal{C}))$ are isomorphic as monoidal triangulated categories.
\end{proposition}
\begin{proof}
	\cite[8.4]{hoveyspec} proved this isomorphism in the case of a weak equivalence $f:T\to S$ in $\mathcal{C}$.
	
	The general case follows by constructing a square
	\[
	\xymatrix{
		T \ar[r] \ar[d] & S \ar[d]
		\\
		T' \ar[r] & S'
	} 
	\]
	of isomorphisms in $Ho\mathcal{C}$, where $T\to T'\to 0$ (resp. $S\to S'\to 0$) is a factorization of $T\to 0$ (resp. $S\to 0$) into trivial cofibration followed by fibration. Since $T$ and $S$ are cofibrant, it follows that $T',S'$ are both fibrant and cofibrant. Thus $Ho\mathcal{C}(T',S')=\mathcal{C}(T',S')/\sim$, so the isomorphism $T'\to S'$ in $Ho\mathcal{C}$ can be lifted to a map in $\mathcal{C}$, which must therefore be a weak equivalence. Now, using the result of the proposition for weak equivalences between cofibrant objects in $\mathcal{C}$, we see that $Ho(Sp_T^\Sigma(\mathcal{C}))\cong Ho(Sp_{T'}^\Sigma(\mathcal{C}))\cong Ho(Sp_{S'}^\Sigma(\mathcal{C}))\cong Ho(Sp_S^\Sigma(\mathcal{C}))$, as desired.
\end{proof} 
\noindent As a more simple observation in the same spirit, given a cofibrant object $S\in \mathcal{C}$ and an integer $i$, the shift $S[i]$ will also be cofibrant and there is an isomorphism of monoidal model categories $Sp^\Sigma_{S}(\mathcal{C})\to Sp^\Sigma_{S[i]}(\mathcal{C})$ mapping the sequence $(E_n)$ to the sequence $(E_n[in])$.
\\
\\Recall that, given a pointed $k$-scheme $(X,x)$, we let $\mathbb{Z}_{tr}(X,x)$ be the direct summand of $\mathbb{Z}_{tr}(X)$ given by the projector $X\to \text{Spec }k \xrightarrow{x} X$. Note that $\mathbb{Z}_{tr}(X,x)$, being a direct summand of the cofibrant object $\mathbb{Z}_{tr}(X)$, is cofibrant as well. We let $G_m$ denote the pointed $k$-scheme $(\mathbb{A}^1-\{0\},1)$, and we let $T$ denote the pointed $k$-scheme $(\mathbb{P}^1,1)$. Then we have an $\mathbb{A}^1$-equivalence $\mathbb{Z}_{tr}(T)\cong \mathbb{Z}_{tr}(G_m)[1]$.
\\
\\Cisinski-Deglise defined $DM(k)$ as the homotopy category of  $Sp_{\mathbb{Z}_{tr}(G_m)}(\mathcal{C})$, but by Proposition 1.1 and the subsequent remark, it follows that the model category $Sp_{Tate}(k):=Sp_{\mathbb{Z}_{tr}(T)}(\mathcal{C})$ has a canonically equivalent homotopy category, so we may work with this model category, which is more convenient for our construction, instead. We may implicitly identify objects of $\mathcal{C}$ with objects of $Sp_{Tate}(k)$ via $\Sigma^\infty$.

Note that $\mathbb{Z}_{tr}(T)\cong \mathbb{Z}(1)[2]$ in $DM(k)$, so we have $R\Omega_T^r(X)\cong X(-r)[-2r]$ in $DM(k)$ for any $X\in Sp_{Tate}(k)$. 

We also recall that for any algebraic space $X$ of finite type over $k$, and smooth $k$-scheme $S$, $z_{equi}(X,r)(S)$ is the group of correspondences from $S$ to $X$ which are equidimensional of relative dimension $r$ over $S$. The $z_{equi}(X,r)$ are Nisnevich sheaves and $z_{equi}(X,0)\cong M^c(X)$ in $DM(k)$.
\section{Functorial construction of  $LM^c:\mathcal{A}\to DM(k)$}
The main difficulty in constructing an analog of the $M^c$ functor in the level of the model category $Sp_{Tate}(k)$ is that given a flat map $f\colon E\to X$ of relative dimension $r$, there is no direct map in $Sp_{Tate}(k)$ which gives the pullback map $M^c(X)\to M^c(E)(-r)[-2r]$ in $DM(k)$. Instead, there is only a zig-zag in $Sp_{Tate}(k)$ inducing the pullback map. For this reason, we must work with zig-zag diagrams as opposed to sequences.
\\
\\Remark: for notational simplicity, we carry out all constructions in the integral case, although everything carries through to a general coefficient ring.
\\
We introduce the following shorthand for the sake of readability of future diagrams. In what follows, let $Q$ denote functorial cofibrant replacement, as in Hovey\cite{hoveybook}.
\begin{definition}
	\label{Mcijr}
	given numbers $i,j,r\in \mathbb{N}$ and algebraic space $X$ of finite type over $k$, let
	\[
	\hat{M}^c_{i,j,r}(X)=(R\Omega_T)^{\circ i} \circ L\Sigma_T^\infty\circ (\Omega_T\circ Q)^{\circ j} (z_{equi}(X,r))
	\]
	In particular, $\hat{M}^c_{i,0,0}(X)$ corresponds in $DM(k)$ to the motive $M^c(X)(-i)[-2i]$.
\end{definition}
\noindent Note: as a convention, we will top our constructions with a hat if they are carried out in the model category rather than the homotopy category.
\begin{lemma}
	\label{mcijrfunctoriality}
	Let $i,j,r,t\in \mathbb{N}$ and let $X, E, Y$ be algebraic spaces of finite type over $k$. Then we have the following canonical maps in $Sp_{Tate}(k)$:
	
	a) $\hat{M}^c_{i,j,r}(X)\to \hat{M}^c_{i,j,r+s}(E)$ for any flat map $E\to X$ of relative dimension $s$.
	
	b) $\hat{M}^c_{i,j+t,r}(X)\to \hat{M}^c_{i,j,r+t}(X)$, which is an $\mathbb{A}^1$-weak equivalence.
	
	c) $\hat{M}^c_{i,t+j,r}(X)\to \hat{M}^c_{i+t,j,r}(X)$.
	
	d) $\hat{M}^c(f):\hat{M}^c_{i,j,r}(X)\to \hat{M}^c_{i,j,r}(Y)$, for any proper map $f:X\to Y$.
\end{lemma}
The maps in a) - d) are the basic maps we will use in our construction.
\begin{proof}
	a) Follows from the natural flat pullback maps $z_{equi}(X,r)\to z_{equi}(E,r+s)$.
	
	b) We can reduce to the case $i,j=0, t=1$. 
	
	Since the cofibrant resolution $Qz_{equi}(X,r)\to z_{equi}(X,r)$ is a trivial fibration, and $\Omega_T$, being right-Quillen, preserves trivial fibrations, it follows that the map $\Omega_T Q z_{equi}(X,r)\to \Omega_T z_{equi}(X,r)$ is an $\mathbb{A}^1$-weak equivalence. Thus it suffices to establish a map $\Omega_T(z_{equi}(X,r))\to z_{equi}(X,r+1)$ which is an $\mathbb{A}^1$-weak equivalence.
	
	Recall that for any algebraic space $X$ of finite type over $k$ and Nisnevich sheaf with transfers $F$, we have $\underline{Hom}_{\mathcal{C}}(\mathbb{Z}_{tr}(X),F)(U)=F(U\times X)$ for all $U\in Sm/k$. Thus there is a natural map $\underline{Hom}_\mathcal{C}(\mathbb{Z}_{tr}(\mathbb{P}^1),z_{equi}(X,r))\to z_{equi}(X,r+1))$ given by
	\begin{equation*}
		\underline{Hom}_\mathcal{C}(\mathbb{Z}_{tr}(\mathbb{P}^1),z_{equi}(X,r))(U)=z_{equi}(X,r)(U\times \mathbb{P}^1)
	\end{equation*}
	\begin{equation*}
		\to z_{equi}(X\times \mathbb{P}^1,r+1)(U)
		\to z_{equi}(X,r+1)(U)
	\end{equation*}
	where the first arrow follows since any cycle $W$ of $U\times \mathbb{P}^1\times X$ which is equidimensional of relative dimension $r$ over $U\times \mathbb{P}^1$ may also be viewed as equidimensional of relative dimension $r+1$ over $U$; and the second arrow follows by proper pushforward. 
	
	Restricting the above map to the direct summand $\mathbb{Z}_{tr}(T)$ of $\mathbb{Z}_{tr}(\mathbb{P}^1)$, and recalling that $\Omega_T=\underline{Hom}_\mathcal{C}(\mathbb{Z}_{tr}(T),-)$, we obtain the desired map 
	\[\Omega_T(z_{equi}(X,r))\to z_{equi}(X,r+1)
	\]
	Voevodsky\cite{voe} proved this is indeed an $\mathbb{A}^1$-weak equivalence in the case $r=0$. The general case follows similarly.
	
	c) It suffices to show that for any object $F\in \mathcal{C}$, there is a natural morphism $L\Sigma_T^\infty(\Omega_T Q F) \to R\Omega_T( L\Sigma_T^\infty F)$. Indeed, to see this we first note that by adjunction, we have a natural map $\Sigma_T^\infty \Omega_T F\to \Omega_T \Sigma_T^\infty F$. Thus we obtain a map
	\[
	L\Sigma_T^\infty(\Omega_TQF)\to \Sigma_T^\infty(\Omega_TQF)\to \Omega_T(\Sigma_T^\infty QF)=\Omega_T(L\Sigma_T^\infty F)\to R\Omega_T(L\Sigma_T^\infty F)
	\] 
	d) By the  proper pushforward $z_{equi}(X,r)\to z_{equi}(Y,r)$.
\end{proof}
Note on functoriality: The maps of Lemma 2.1a)-d) behave well with each other in the following ways. The maps of type a) are functorial with respect to flat maps, and the maps of type d) are functorial with respect to proper maps. The maps of type b) or c) commute with any of the maps of type a)-d). And finally, the pullback maps a) commute with the pushforward maps d) for any fiber square
\[
\xymatrix{
	E_1 \ar[d] \ar[r] & E_2 \ar[d]\\
	X_1 \ar[r] & X_2
}
\]
where the horizontal maps are proper and the vertical maps are flat.
\\
\\
The following is the model-category version of flat pullback, which forms the main ingredient of our construction.
\begin{construction}
	\label{mcijrf}
	Let $i,j,r\in \mathbb{N}$. Let $f:E\to X$ be a flat map of algebraic spaces of finite type over $k$, of relative dimension $d$. Then we define $\hat{M}^c_{i,j,r}(f)$ as the following zig-zag in $Sp_{Tate}(k)$:
	\[
	\xymatrix{
		M^c_{i,j,r}(X) \ar[dr]^{\alpha} & & M^c_{i,j+d,r}(E) \ar[dl]_{\beta} \ar[dr]^{\gamma} &\\
		& M^c_{i,j,r+d}(E) & & M^c_{i+d,j,r}(E)
	}
	\]
	Here $\alpha,\beta,\gamma$ are given by Lemma 2.1 a),b),c), respectively.
\end{construction}
\noindent Note: By Lemma \ref{mcijrfunctoriality}, the map $\beta$ is an $\mathbb{A}^1$-weak equivalence. So the zig-zag $\hat{M}^c_{i,j,r}(f)$ represents a map $\gamma \circ \beta^{-1}\circ \alpha$ in $DM(k)$ from $M^c_{i,j,r}(X)$ to $M^c_{i+t,j,r}(E)$. We denote this map $M^c_{i,j,r}(f)$. 

\begin{lemma}
	\label{pullback}
	a) As a special case, $M^c_{i,0,0}(f)$ is the natural flat pullback map $M^c(X)(-i)[-2i]$ to $M^c(E)(-i-t)[-2i-2t]$.
	
	b) If furthermore $f$ is a vector bundle, then $\alpha$ and $\gamma$ in the construction of $\hat{M}^c_{i,0,0}(f)$ are $\mathbb{A}^1$-weak equivalences.
\end{lemma}
\noindent Note: it follows that the reverse zig-zag of $\hat{M}^c_{i,0,0}(f)$, which we will denote $\hat{M}^c_{i,0,0}(f)^{-1}$, corresponds in $DM(k)$ to the isomorphism $M^c(E)(-i-t)[-2i-2t]\to M^c(X)(-i)[-2i]$ inverse to the natural pullback.
\begin{proof} a) Follows from \cite{voe}.
	%todo: do this part better
	
	b) We first show that $\gamma$ is an $\mathbb{A}^1$- weak equivalence, for which it suffices to show that the map $M^c_{0,d,0}(E)\to M^c_{d,0,0}(E)$ is. Consider the following commutative diagram in $DM(k)$:
	
	\[
	\xymatrix{
		M^c_{0,d,0}(E)=L\Sigma_T^\infty((\Omega_T\circ Q)^{\circ d}(z_{equi}(E,0))) \ar[r] \ar[d]^{\epsilon} & M^c_{d,0,0}(E)=(R\Omega_T)^{\circ d}(L\Sigma_T^\infty z_{equi}(E,0)) \ar@{=}[dd]\\
		L\Sigma_T^\infty (R\Omega_T^d(z_{equi}(E,0))) \ar@{=}[d] & \\
		\underline{Hom}_{DM^{eff}(k)}(\mathbb{Z}(d)[2d],M^c(E)) \ar[r]_{\theta} & \underline{Hom}_{DM(k)}(\mathbb{Z}(d)[2d],M^c(E))
	}
	\]
	
	where all maps are the natural ones (and $\theta$ is only defined in $DM(k)$). 
	
	Voevodsky\cite{voe} showed that the natural map $\Omega_T^dz_{equi}(E,0)\to R\Omega_T^dz_{equi}(E,0)$ is an isomorphism in $DM^{eff}(k)$, hence so is $\epsilon$. 
	
	On the other hand, because $E\to X$ is a vector bundle it follows that we have $\underline{Hom}_{DM(k)}(\mathbb{Z}(d)[2d],M^c(E))\cong M^c(E)(-d)[2d]\cong M^c(X)\in DM^{eff}(k)$, so in fact $\theta$ is an isomorphism. Thus $\alpha$, the diagonal of the above square, induces an isomorphism in $DM(k)$, so it is an $\mathbb{A}^1$-weak equivalence.
	\\
	\\Now since $\beta$ and $\gamma$ are both isomorphisms in $DM(k)$ and $\gamma \circ \beta^{-1}\circ \alpha=M^c(f)$ is an isomorphism in $DM(k)$ as $f$ is a vector bundle, it follows that $\alpha$ must be an isomorphism in $DM(k)$, hence it is an $\mathbb{A}^1$-weak equivalence in $\mathcal{C}$. \end{proof}
\noindent If $X$ is an algebraic space of finite type over $k$ and $\iota:U\to X$ is an open embedding, then $\iota$ is a flat map of relative dimension $0$. But in this case the pullback $M^c_{i,j,r}(\iota)$ can be represented in $Sp_{Tate}(k)$ directly as a morphism $\hat{M}^c_{i,j,r}(\iota):\hat{M}^c_{i,j,r}(X)\to \hat{M}^c_{i,j,r}(U)$, as opposed to a zig-zag.
\\
\\Now we are ready to extend the $\hat{M}^c_{i,j,r}$ construction to a similar construction $\hat{M}^c_{i,j,r}:\mathcal{A}\to Sp_{Tate}(k)^{ZigZag}$, where $ZigZag$ is the indexing category
\[
\xymatrix@R=0.5em{
	\cdot \ar[ddr]  & & \cdot \ar[ddl] \ar[ddr]& & \cdot \ar[ddl]
	\\
	\cdots & & & &
	\\
	& \cdot & & \cdot &
}
\]
\begin{definition}
	\label{mcijrhat}
	given an object $X\in \mathcal{A}$ of the form 
	\[
	\cdots \to X_3 \xleftarrow{i_2} Y_2\xrightarrow{j_2} E_2\xrightarrow{f_2} X_2 \xleftarrow{i_1} Y_1\xrightarrow{j_1} E_1\xrightarrow{f_1} X_1
	\]
	where $f_i$ is flat of relative dimension $d_i$, define $\hat{M}^c_{i,j,r}(X)$ to be the concatenation (denoted by $\diamond$) of zig-zags
	\[
	\cdots \diamond \hat{M}^c_{i+d_1,j,r}(f_2)^{-1}\diamond (\hat{M}^c_{i+d_1,j,r}(j_1)\circ \hat{M}^c_{i+d_1,j,r}(i_1))\diamond \hat{M}^{c}_{i,j,r}(f_1)^{-1}
	\]
	outlined as follows:
	\[
	\xymatrix@R=0.5em@C=0.5em{
		& \cdot \ar[ddr] \ar[ddl] & & M^c_{i+d_1,j,r}(X_2) \ar[ddr] \ar[ddl]  & &  \ar[ddl] \cdot \ar[ddr]& & M^c_{i,j,r}(X_1) \ar[ddl] &
		\\
		\cdots & & & & & & & &
		\\
		M^c_{i+d_1+d_2,j,r}(E_2) & & \cdot & & M^c_{i+d_1,j,r}(E_1) & & \cdot & & 
	}
	\]
\end{definition} 
\noindent Note: Define $\hat{M}^c:\mathcal{A}\to Sp_{Tate}(k)^{ZigZag}$ to be $\hat{M}^c_{0,0,0}$. Then by lemma 2.3, it induces the functor $M^c:\mathcal{A}\to DM(k)^{\mathbb{N}^{op}}$ upon localization, inverting the $\mathbb{A}^1$-weak equivalences and passing to a subsequence.
\\
\\We know that $M^c$ is functorial on $\mathcal{A}$. To study the functoriality of $\hat{M}^c$, we simply iterate our above constructions. The following lemma is an analogue of Lemma \ref{mcijrfunctoriality} in the context of $\hat{M}^c$.
\begin{lemma}
	\label{mcijrhatfunctoriality}
	Let $i,j,r,t\in \mathbb{N}$ and let $X, E, Y\in \mathcal{A}$. Then we have the following canonical maps in $Sp_{Tate}(k)^{ZigZag}$:
	
	a) $\hat{M}^c_{i,j,r}(X)\to \hat{M}^c_{i,j,r+s}(E)$ for any flat map $E\to X$ of relative dimension $s$.
	
	b) $\hat{M}^c_{i,j+t,r}(X)\to \hat{M}^c_{i,j,r+t}(X)$, which is a termwise $\mathbb{A}^1$-weak equivalence.
	
	c) $\hat{M}^c_{i,t+j,r}(X)\to \hat{M}^c_{i+t,j,r}(X)$.
	
	d) $\hat{M}^c(f):\hat{M}^c_{i,j,r}(X)\to \hat{M}^c_{i,j,r}(Y)$, for any proper map $f:X\to Y$.
\end{lemma}
\begin{proof}
	We apply the corresponding maps from Lemma 2.1 termwise to the relevant zig-zags. These maps form a morphism of diagrams by the functoriality properties of Lemma 2.1.
\end{proof}
\noindent Note: The maps in Lemma 2.5 have the same functoriality properties as noted for Lemma 2.1.
\\
\\We can similarly iterate Construction \ref{mcijrf}:
\begin{construction}
	Let $f:E\to X$ be flat of relative dimension $d$ in $\mathcal{A}$. 
	Then we define $\hat{M}^c_{i,j,r}(f)$ to be the following zig-zag in $Sp_{Tate}(k)^{ZigZag}$
	\[
	\xymatrix{
		\hat{M}^c_{i,j,r}(X) \ar[dr]^{\alpha} & & \hat{M}^c_{i,j+d,r}(E) \ar[dl]_{\beta} \ar[dr]^{\gamma} &\\
		& \hat{M}^c_{i,j,r+d}(E) & & \hat{M}^c_{i+d,j,r}(E)
	}
	\]
	where the morphisms are given by Lemma \ref{mcijrhatfunctoriality}
\end{construction}
\noindent Having defined $\hat{M}^c:\mathcal{A}\to Sp_{Tate}(k)^{ZigZag}$, we can now take the homotopy limit in the setting of a model category which, unlike with triangulated categories, is functorial.

It helps to consider $Sp_{Tate}(k)^{ZigZag}$ as a model category. For any model category $\mathcal{C}$ and an inverse small category $D$, the diagram category $\mathcal{C}^D$ has the Reedy model category structure. A weak equivalence in $Sp_{Tate}(k)^{ZigZag}$ is a termwise weak-equivalence. For our purposes we only need to additionaly use the fact that a fibrant object in $Sp_{Tate}(k)^{ZigZag}$ is a zig-zag consisting of fibrations between fibrant objects. The limit functor $lim\colon Sp_{Tate}(k)^{ZigZag}\to Sp_{Tate}(k)$ is a right-Quillen functor. We define $holim\colon Ho(Sp_{Tate}(k)^{ZigZag})\to DM(k)$ to be its right-derived functor.
\\
\\Given an object $X\in \mathcal{A}$, we let $\tilde{M}^c_{i,j,r}(X)$ denote the image of $\hat{M}^c_{i,j,r}(X)$ in the localization $Ho(Sp_{Tate}(k)^{ZigZag})$. We will focus in particular on the case $\tilde{M}^c_{i,0,0}(X)$, which we may abbreviate $\tilde{M}^c_i(X)$. Let $LM^c_i(X)\in DM(k)$ denote $holim(\tilde(M^c_i(X)))$.

Let $f:E\to X$ be a flat map in $\mathcal{A}$ of relative dimension $d$. Since the map $\beta$ in Construction 2.6 is termwise an $\mathbb{A}^1$-weak equivalence, the zig-zag $\hat{M}^c_{i,0,0}(f)$ above induces a map $\tilde{M}^c_{i}(f): \tilde{M}^c_{i}(X)\to \tilde{M}^c_{i+d}(E)$ given by the composition $\gamma\circ\beta^{-1}\circ\alpha$. Let $LM^c_i(f):LM^c_i(X)\to LM^c_{i+d}(E)$ be the result of applying the $holim$ functor to $\tilde{M}^c_i(f)$.

Similarly, if $f:X\to Y$ is a proper map in $\mathcal{A}$, let $\tilde{M}^c_{i}(f):\tilde{M}^c_{i}(X)\to \tilde{M}^c_{i}(Y)$ be the localization of the map $\hat{M}^c_{i,0,0}(f): \hat{M}^c_{i,0,0}(X)\to \hat{M}^c_{i,0,0}(Y)$. Let $LM^c_i(f):LM^c_i(X)\to LM^c_i(Y)$ be the result of applying the $holim$ functor to $\tilde{M}^c_i(f)$.
\begin{proposition}
	\label{limitfunctoriality}
	(Functoriality of the $LM^c_i$): 
	
	a)(Functoriality with respect to proper maps). Let $f:X\to Y$ and $g:Y\to Z$ be proper maps in $\mathcal{A}$. Then $LM^c_i(g)\circ LM^c_i(f)=LM^c_i(g\circ f)$ for all $i\in \mathbb{N}$.
	
	b)(Functoriality with respect to flat maps). Let $\theta:F\to E$ and $\varphi:E\to X$ flat maps in $\mathcal{A}$ of relative dimensions $e$ and $d$, respectively. Then $LM^c_i(\varphi \circ \theta)=LM^c_{i+d}(\theta)\circ LM^c_i(\varphi)$ for all $i\in \mathbb{N}$
	
	c)(Mixing flat and proper maps) Consider a square
	\[
	\xymatrix{
		E_1 \ar[r]^{f} \ar[d]^{\varphi} & E_2 \ar[d]^{\theta}\\
		X_1 \ar[r]^{g} & X_2
	}
	\]
	in $\mathcal{A}$ where the horizontal maps are proper, the vertical maps are flat of relative dimension $d$, and the square is componentwise a fiber square. 
	
	Then $ LM^c_i(\theta)\circ LM^c_i(g)=LM^c_{i+d}\circ LM^c_i(\varphi)$
\end{proposition}
\begin{proof}
	We prove the functoriality properties in a)-c) for the $\tilde{M}^c_i$ functors. They will then follow for the $LM^c_i$ functors by applying the $holim$ functor. We also reduce to the case $i=0$ throughout.
	\\
	\\a) Follows from the functoriality of $\hat{M}^c_{i,j,r}$ with respect to proper maps.
	\\
	\\b)Consider the following diagram in $Sp_{Tate}(k)^{ZigZag}$, where all morphisms are the natural ones given by Lemma 2.5 and applying the holim functor:
	\[
	\xymatrix@C=0.3em@R=0.8em{
		\hat{M}^c_{0,0,0}(X) \ar[ddr] \ar[dddddr] & & \hat{M}^c_{0,d,0}(E) \ar[ddd] \ar[ddl] \ar[ddr] & & \hat{M}^c_{d,e,0}(F) \ar[dddddl] \ar[ddr] &\\
		& & & & &\\
		& \hat{M}^c_{0,0,d}(E) \ar[ddd] & & \hat{M}^c_{d,0,0}(E) \ar[ddd] & & \hat{M}^c_{d+e,0,0}(F)
		\\& & \hat{M}^c_{0,d,e}(F) \ar[ddl] \ar[ddr] & & &
		\\& & & & &
		\\& \hat{M}^c_{0,0,d+e}(F) & & \hat{M}^c_{d,0,e}(F) & &
	}
	\]
	The commutativity of the diagram demonstrates that the upper zig-zag and the lower zig-zag from $\hat{M}^c_{0,0,0}(X)$ to $\hat{M}^c_{d+e,0,0}(F)$ induce the same map in the homotopy category $Ho(Sp_{Tate}(k)^{ZigZag})$.
	
	The upper zig-zag  is the concatenation $\hat{M}^c_{0,0,0}(\varphi) \diamond \hat{M}^c_{d,0,0}(\theta)$. On the other hand, the lower zig-zag itself induces the same map as $\hat{M}^c_{0,0,0}(\varphi\circ \theta)$ in $Ho(Sp_{Tate}(k)^{ZigZag})$, by similarly comparing the top zig-zag and bottom zig-zag in the following commutative diagram:
	\[
	\xymatrix@C=0.3em@R=0.8em{
		& & \hat{M}^c_{0,d+e,0}(F) \ar[ddddl] \ar[dd] \ar[ddrr]
		\\
		\\
		\hat{M}^c_{0,0,0}(X) \ar[ddr] & & \hat{M}^c_{0,d,e}(F) \ar[ddl] \ar[ddr] & & \hat{M}^c_{d,e,0}(F) \ar[ddl] \ar[ddr] &\\
		& & & & &\\
		& \hat{M}^c_{0,0,d+e}(F) & & \hat{M}^c_{d,0,e}(F) & & \hat{M}^c_{d+e,0,0}(F)
	}
	\]
	In conclusion, the concatenation $\hat{M}^c_{0,0,0}(\varphi) \diamond \hat{M}^c_{d,0,0}(\theta)$ induces the same map as $\hat{M}^c_{0,0,0}(\varphi\circ \theta)$ in $Ho(Sp_{Tate}(k)^{ZigZag})$, which is to say  that indeed $\tilde{M}^c_{d,0,0}(\theta)\circ \tilde{M}^c_{0,0,0}(\varphi) =\tilde{M}^c_{0,0,0}(\varphi\circ \theta)$, as desired. 
	\\
	\\c) Consider the natural commutative diagram
	\[
	\xymatrix{
		\hat{M}^c_{0,0,0}(X_1) \ar[rd] \ar[dd] & & \hat{M}^c_{0,d,0}(E_1) \ar[dl] \ar[dr] \ar[dd] &\\
		& \hat{M}^c_{0,0,d}(E_1) \ar[dd] & & \hat{M}^c_{d,0,0}(E_1) \ar[dd]\\
		\hat{M}^c_{0,0,0}(X_2) \ar[rd] & & \hat{M}^c_{0,d,0}(E_2) \ar[dl] \ar[dr] &\\
		& \hat{M}^c_{0,0,d}(E_2)  & & \hat{M}^c_{d,0,0}(E_2)
	}
	\]
	By commutativity of the diagram, the top path and the bottom path from $\hat{M}^c_{0,0,0}(X_1)$ to $\hat{M}^c_{d,0,0}(E_2)$ induce the same map in $Ho(Sp_{Tate}(k)^{ZigZag})$. 
	
	That is, $\tilde{M}^c_d(f)\circ \tilde{M}^c_0(\varphi)=\tilde{M}^c_0(\theta)\circ\tilde{M}^c_0(g)$, as desired.
\end{proof}
\noindent Having shown the functoriality of the $LM^c_i$, it remains to show that $LM^c_i$, which we defined using a homotopy limit on the model category level, agrees with the previously defined $LM^c$ construction, which was defined using a homotopy limit on the triangulated category level.
\begin{proposition}
	\label{limitsagree}
	Let $X\in \mathcal{A}$, $i\in \mathbb{N}$. Then $LM^c_{i,0,0}(X)=LM^c(X)(-i)[-2i]$, where $LM^c(X)$ is as defined in \ref{McA}
\end{proposition}
\noindent To prove Proposition \ref{limitsagree}, we examine the relationship between the homotopy limit constructions in model categories and the analogous one in triangulated categories. Working in generality and overriding previous notation, let $\mathcal{C}$ denote an additive model category.

We begin by comparing the additive structure with the model category structure.
\begin{lemma}
	\label{additivemodel}
	a) For $\mathcal{W}$ equal to any of the classes fibrations, cofibrations, trivial fibrations, trivial cofibrations, or weak equivalences in $\mathcal{C}$, the class $\mathcal{W}$ of maps is closed under direct sum, in the sense that if $f:A\to B$ and $g:C\to D$ are in $\mathcal{W}$, then so is $f\oplus g:A\oplus C\to B\oplus D$. 
	\\b) The localization functor $\gamma:\mathcal{C}\to Ho(\mathcal{C})$ is additive.
\end{lemma}
\begin{proof}
	a) In a general model category, fibrations and trivial fibrations are closed under coproduct. Similarly, in a general model category cofibrations and trivial cofibrations are closed under product. Thus in an additive category, all of these classes of maps are closed under direct sum. It then follows that weak equivalences are closed under direct sum, since every weak equivalence can be factored into a trivial fibration and a trivial cofibration.
	\\
	\\b) In a general model category, given an arbitrary index set $I$, the product map $\prod_I:\mathcal{C}^I\to \mathcal{C}$ is a right Quillen functor adjoint to the diagonal map, giving rise to a right-derived functor on the homotopy categories, $R\prod_I: Ho\mathcal{C}^I\to Ho\mathcal{C}$, adjoint to the diagonal map. This adjunction shows that $R\coprod_I$ is the product in $Ho\mathcal{C}$. Thus if $A_i: i\in I$ are fibrant objects in $\mathcal{C}$, then $\gamma(\prod_{i\in I} A_i) = \prod_{i\in I}\gamma(A_i)$.
	
	Dually, the left-derived functor $L\coprod_I: Ho\mathcal{C}^I\to Ho\mathcal{C}$ is the product in $Ho\mathcal{C}$, so if $A_i: i\in I$ are cofibrant objects in $\mathcal{C}$, then $\gamma(\coprod_{i\in I} A_i) = \coprod_{i\in I}\gamma(A_i)$. 
	
	Now suppose $\mathcal{C}$ is additive and let $A,B\in \mathcal{C}$. Since weak equivalences are closed under direct sum by a), we see by replacing $A$ and $B$ by fibrant objects that $\gamma(A\oplus B)=\gamma(A)\coprod \gamma(B)$. On the other hand, replacing $A$ and $B$ by cofibrant objects shows that $\gamma(A\oplus B)=\gamma(A)\prod \gamma(B)$. Thus, the biproduct of $\gamma(A)$ and $\gamma(B)$ exists and equals $\gamma(A\oplus B)$. By similarly considering the empty biproduct, we can see that $\gamma(0)=0$. 
	
	Thus $Ho\mathcal{C}$ is at least semi-additive, so Hom groups $Ho \mathcal{C}(\gamma A,\gamma B)$ are at least abelian monoids. To show additive inverses, we may assume by fribrant and cofibrant replacement that $A$ is fibrant and $B$ cofibrant. In this case we have that $\gamma: \mathcal{C}(A,B)\to Ho\mathcal{C}(\gamma A, \gamma B)$ is quotient by homotopy, so it is surjective. Since $\gamma$ is a homomorphism of monoids and $\mathcal{C}(A,B)$ has additive inverses, the monoid $Ho \mathcal{C}(\gamma A, \gamma B)$ must have additive inverses as well. So $Ho \mathcal{C}$ is additive and $\gamma$ is an additive functor.
\end{proof}
\begin{lemma}
	\label{shiftfibrations}
	Let $\cdots A_2\to A_1\to A_0\to 0$ be a sequence of fibrations. Then the map $\text{id}-\text{shift}:\prod_i A_i \to \prod_i A_i$ is a fibration as well.
\end{lemma}
\begin{proof} We prove this by verifying that $\text{id}-\text{shift}$ satisfies the right-lifting property with respect to trivial cofibrations. Consider a diagram
	\[
	\xymatrix{
		X \ar[r]^{(t_i)} \ar[d]^{s} & \prod_i A_i \ar[d]^{\text{id}-\text{shift}} \\
		Y \ar[r]^{(r_i)} & \prod_i A_i
	}
	\]
	where $s:X\to Y$ is a trivial cofibration. Then we recursively construct the following commutative diagrams:
	\[
	\xymatrix{
		X \ar[r]^{t_0} \ar[d]^{s} & A_0
		\\
		Y\ar[ur]_{h_1}  
	},
	\xymatrix{
		X \ar[r]^{t_{i+1}} \ar[d]^{s} & A_{i+1} \ar[d]^{f_{i+1}}
		\\
		Y\ar[ur]_{h_{i+1}} \ar[r]_{h_i-r_i} & A_i  
	}
	\]
	Where the lifts $h_i$ exist by the lifting properties of the $f_i$. Then $(h_i):Y\to \prod_i A_i$ is a lift of the original diagram, since for each $i$, we have $h_i \circ s=t_i$ and $h_i-f_{i+1}\circ h_{i+1}=r_i$.
\end{proof}
\noindent Now assume further that $\mathcal{C}$ is stable, so the homotopy category $Ho(\mathcal{C})$ is triangulated. Given a sequence of maps $\cdots \to A_2\to A_1\to A_0$ in a triangulated category $\mathcal{D}$ which has arbitrary products, the homotopy limit $\text{holim}A_i$ in the sense of Bokstedt-Neeman\cite{bokneem} is defined as the fiber of the morphism $\text{id}-\text{shift}:\prod_i A_i\to \prod_i A_i$
\begin{lemma}
	\label{limitfibrations}
	a) Given a diagram $\cdots \to A_2\to A_1\to A_0\to 0$ in $\mathcal{C}$ consisting of fibrations, then the limit $\text{lim}A_i$, taken in $\mathcal{C}$, is isomorphic to the homotopy limit $\text{holim}A_i$ taken in $Ho \mathcal{C}$.
	\\
	b) Let $F$ be a diagram $\mathcal{C}^{ZigZag}$ of the form
	\[
	\xymatrix@R=0.5em{
		A_3 \ar[ddr]^{s_2} & & A_2 \ar[ddl]^{t_2} \ar[ddr]_{s_1}& & A_1 \ar[ddl]^{t_1}
		\\
		\cdots & & & &
		\\
		& B_2 & & B_1 & 
	}
	\]
	where all maps are fibrations between fibrant objects and the $t_i$ are trivial fibrations. Then \[\text{lim}F=\text{holim}(\cdots\to A_3\xrightarrow{ t_2^{-1}\circ s_2}A_2\xrightarrow{t_1^{-1}\circ s_1}A_1)\].
\end{lemma}
\begin{proof}
	a) As in the proof of lemma \ref{additivemodel}b), since the $A_i$ are all fibrant, then the product $\prod_i A_i$, taken in $\mathcal{C}$, is still fibrant, and it equals the product $\prod_i A_i$ taken in $Ho\mathcal{C}$. Also, since the localization functor $\mathcal{C}\to Ho\mathcal{C}$ is additive, $\text{id}-\text{shift}:\prod_i A_i\to \prod_i A_i$ maps to $\text{id}-\text{shift}$ under the localization. By lemma \ref{shiftfibrations}, the map $\text{id}-\text{shift}$ is a fibration. Thus we obtain a fiber sequence $\text{ker}(\text{id}-\text{shift})\to \prod_i A_i\xrightarrow{\text{id}-\text{shift}} \prod_i A_i$. But $\text{ker}(\text{id}-\text{shift})=\text{lim}A_i$ so we have a fiber sequence $\text{lim}A_i\to \prod_i A_i\xrightarrow{\text{id}-\text{shift}} \prod_i A_i$. Since fiber sequences in stable model categories induce exact triangles in the homotopy category, we see that indeed $\text{lim}A_i\cong \text{holim}A_i$ in $Ho\mathcal{C}$, as desired.
	\\
	\\b) By taking successive pullbacks, form the following diagram:
	\[
	\xymatrix{
		& & X_3 \ar[rr]^{f_3} \ar[d]^{h_3} & & X_2 \ar[rr]^{f_2} \ar[d]^{h_2}& & X_1\ar@{=}[d]^{h_1} 
		\\
		\cdots  & & A_3 \ar[dl]^{t_3} \ar[dr]_{s_2} & & A_2 \ar[dl]^{t_2} \ar[dr]_{s_1}& & A_1 \ar[dl]^{t_1} 
		\\
		& B_3 & & B_2 & & B_1 &
	}
	\]
	That is, define $X_1=A_1$, $X_{n+1}=A_{n+1}\times_{B_{n}} X_n$.
	\\
	\\Since fibrations and trivial fibrations are closed under pullback and composition, it follows that the $X_i$ are fibrant, the maps $f_{i+1}:X_{i+1}\to X_i$ are fibrations, and the maps $h_i:X_i\to A_i$ are trivial fibrations, thus forming an isomorphism between the sequences $\cdots X_3\to X_2\to X_1$ and $\cdots A_3\to A_2\to A_1$ in $Ho\mathcal{C}$.
	
	By part a), $\text{lim}X_i \cong \text{holim}X_i$ in $Ho\mathcal{C}$. But $\text{lim}F=\text{lim}X_i$ by comparing the two limits' universal properties, and $(X_i)\cong (A_i)$ in $Ho\mathcal{C}$. Thus  $\text{lim}F=\text{holim}A_i$, as desired.
\end{proof}
\noindent Proof of Proposition \ref{limitsagree}: Let $X\in \mathcal{A}$. Applying Proposition \ref{limitfibrations} to the category $\mathcal{C}=Sp_{Tate}(k)$ and the zig-zag $F=M^c_{i,0,0}(X)$, we obtain Proposition \ref{limitsagree} as a corollary. Indeed, keeping the notation of Proposition \ref{limitfunctoriality}b) applied to this case, we use that the sequence $\cdots\to A_5\to A_3\to A_1$ is precisely the sequence $M^c(X)(-i)[-2i]$, and homotopy limit is preserved by passing to a subsequence. \qed
\section{Localization triangles}
Define an open embedding in $\mathcal{A}$ to be a flat morphism $i:U\to X$ in $\mathcal{A}$ which is termwise an open embedding. Note that by our definition of flat morphisms, $i$ is of the following form:
\[
\xymatrix{
	\cdots \ar[r] & U_2 \ar[d] & Y_1\cap E_1|_{U_1}  \ar[d] \ar[l] \ar[r] & E_1|_{U_1} \ar[r] \ar[d] & U_1 \ar[d] \\
	\cdots \ar[r] & X_2 & Y_1 \ar[l] \ar[r] & E_1 \ar[r] & X_1
}
\]
Similarly, we define a closed embedding in $\mathcal{A}$ to be a proper morphism $j:Z\to X$ in $\mathcal{A}$ which is termwise a closed embedding. By our definition of proper morphisms, $j$ is of the following form:
\[
\xymatrix{
	\cdots \ar[r] & Z_2 \ar[d] & Z_2\cap Y_1  \ar[d] \ar[l] \ar[r] & E_1|_{Z_1} \ar[r] \ar[d] & Z_1 \ar[d] \\
	\cdots \ar[r] & X_2 & Y_1 \ar[l] \ar[r] & E_1 \ar[r] & X_1
}
\]
We now define the closed complement of $i:U\to X$ as above to be the following closed embedding, denoted $j:X-U\to X$:
\[
\xymatrix{
	\cdots \ar[r] & X_2-U_2 \ar[d] & Y_1-Y_1\cap U_2  \ar[d] \ar[l] \ar[r] & E_1|_{X_1-U_1} \ar[r] \ar[d] & X_1-U_1 \ar[d] \\
	\cdots \ar[r] & X_2 & Y_1 \ar[l] \ar[r] & E_1 \ar[r] & X_1
}
\]
Note that $Y_1\cap E_1|_{U_1}\subset U_2$, so $Y_1-Y_1\cap U_2 \subset Y_1-Y_1\cap E_1|_{U_1}$. Thus indeed $Y_1-Y_1\cap U_2\subset E_1|_{X_1-U_1}$, so that the diagram makes sense.
\begin{lemma}
	\label{holimtriangles}
	Let $A\to B\to C$ be a complex of objects in $Sp_{Tate}(k)^{ZigZag}$ which termwise extend to exact triangles in $DM(k)$. Then the homotopy limits of these diagrams also extend to an exact triangle
	\[
	holim(A)\to holim(B)\to holim(C)\to holim(A)[1]
	\]
	in $DM(k)$
\end{lemma}
\begin{proof}
	Since $holim\colon Ho(Sp_{Tate}(k)^{ZigZag})\to DM(k)$ is exact, the proposition will follow once we show that the sequence $A\to B\to C$ is an exact triangle in $Ho(Sp_{Tate}(k)^{ZigZag})$.
	
	We use the description of exact triangles as cofiber sequences. Fix a square
	\[
	\xymatrix{
		W_A \ar[r]^{f} \ar[d]& W_B\ar[d]\\
		A \ar[r] & B
	}
	\]
	in $Sp_{Tate}(k)^{ZigZag}$, where $f:W_A\to W_B$ is a cofibration between cofibrant objects, and the vertical maps are $\mathbb{A}^1$-weak equivalences. Let $P=coker(f)$. Then by definition, $W_A\to W_B\to P$ is a cofiber sequence.
	
	By the assumption that $A\to B\to C$ is a complex, we can extend the above square to a commutative diagram
	\[
	\xymatrix{
		W_A \ar[r]^{f} \ar[d]& W_B \ar[d] \ar[r] & P \ar[d]^{\varphi}\\
		A \ar[r] & B \ar[r] & C
	}
	\]
	But now we note that the upper sequence is termwise a cofiber sequence since $f$, being a cofibration in $Sp_{Tate}(k)^{ZigZag}$, is termwise a cofibration; and the bottom sequence is also termwise an exact sequence in $DM(k)$ by assumption. Moreover, the two leftmost vertical maps are termwise isomorphisms in $DM(k)$. Thus there is a termwise isomorphism $\psi:P_\alpha\to C_\alpha$ in $DM(k)$ for any $\alpha\in ZigZag$. To see that $\varphi$ itself gives such an isomorphism, consider the diagram in $DM(k)$
	\[
	\xymatrix{
		(W_A)\alpha \ar[r]^{f_\alpha} \ar[d] & (W_B)_\alpha \ar[r] \ar[d] & P_\alpha \ar[d]^{\psi}\\
		RA_\alpha \ar[r] & RB_\alpha \ar[r] & RC_\alpha
	}
	\]
	where $R$ denotes fibrant replacement and $\psi$ is an isomorphism. But the top row consists of cofibrant objects and the bottom row of fibrant ones. Thus we can regard $\psi:P_\alpha\to RC_\alpha$ as a morphism in $Sp_{Tate}(k)$, such that the diagram commutes in $Sp_{Tate}(k)$. But since $P_\alpha=\text{coker}(f_\alpha)$, there is a unique morphism  $P_\alpha\to C_\alpha$ making the diagram commute. This implies that $\psi$ equals $\varphi_\alpha$ followed by fibrant replacement, so $\varphi$ is indeed termwise an isomorphism in $DM(k)$, hence a weak equivalence in $Sp_{Tate}^{ZigZag}$.
	
	Thus the sequence $A\to B\to C$ is isomorphic in $Ho(Sp_{Tate}(k)^{ZigZag})$ to the exact triangle $W_A\to W_B\to P$, so it is exact, as desired.
\end{proof}
\begin{corollary}
	\label{localizationstacks}
	Let $i:U\to X$ be an open embedding of quotient stacks with closed complement $j:Z\to X$. Then there is a natural exact triangle in $DM(k)$ of the form
	\[
	M^c(Z)\to M^c(X)\to M^c(U)
	\]
\end{corollary}
\begin{proof}
	Fix a representative $(V_i-S_i)$ is of $X$ in $\mathcal{A}$. Then 
	\[ LM^c((V_i-S_i)|_Z)\to LM^c((V_i-S_i))\to LM^c((V_i-S_i)|_U)
	\]
	is exact since it termwise consists of localization exact triangles.
\end{proof}
\noindent Note: the following section will show that the compactly supported motive $M^c(X)$ does not depend (up to canonical isomorphism) on the choice of presentation $(V_i-S_i)$ of $X$ by algebraic spaces.
\\
\\
The following argument is made by Totaro\cite{burt}, but we repeat it for clarity and to emphasize the canonical nature of the stated isomorphism.
\begin{corollary}
	\label{vanishingcodim}
	Let $i:U\to X$ be an open embedding in $\mathcal{A}$ of the form
	\[
	\xymatrix{
		\cdots \ar[r] & U_2 \ar[d] & Y_1\cap E_1|_{U_1}  \ar[d] \ar[l] \ar[r] & E_1|_{U_1} \ar[r] \ar[d] & U_1 \ar[d] \\
		\cdots \ar[r] & X_2 & Y_1 \ar[l] \ar[r] & E_1 \ar[r] & X_1
	}
	\]
	such that the codimension of $(X_i-U_i)$ in $X_i$ approaches infinity. Then $i$ induces by flat pullback an isomorphism $LM^c(X)\cong LM^c(U)$.
\end{corollary}
\begin{proof}
	We claim that $LM^c(X-U)\in D_j$ for all integers $j$, where $D_j$, as defined by Totaro\cite{burt}, is the smallest localizing category of $DM(k)$ containing $M^c(Y)(a)$ for all $k$-schemes $Y$ and all integers $a$ such that $dim Y+a \le j$. This would imply, as explained by Totaro\cite{burt}, that $LM^c(X-U)=0$, so that by the localization exact triangle, $LM^c(X)\cong LM^c(U)$. 
	
	To prove the claim, let $Z_i:=(X_i-U_i)$. By definition, $LM^c(X-U)$ is the homotopy limit $holim(M^c(Z_i)(-d_1-\cdots-d_{i-1})[-2d_1-\cdots-2d_{i-1}])$, where $d_i$ is the rank of the vector bundle $E_i$ over $X_i$. Here we're viewing it as a homotopy limit of a sequence in $DM(k)$ rather than in the model category. 
	
	Let $e_i=dim(Z_i)-d_1-\cdots-d_{i-1}$. We have \[e_i=dim(X_i)-d_1-\cdots-d_{i-1}-codim(Z_i\subset X_i)\]. 
	But $dim(X_i)-d_1-\cdots-d_{i-1}\le dimX_1<\infty$. So, since $codim(Z_i\subset X_i)$ approaches infinity, we see that $e_i\to -\infty$. By \ref{filtration}, it follows that for every $j\in \mathbb{Z}$, the sequence $(M^c(Z_i)(-d_1-\cdots-d_{i-1})[-2d_1-\cdots-2d_{i-1}])$ is eventually in $D_j$, so indeed its homotopy limit $LM^c(X-U)$ is in $D_j$.
\end{proof}
\section{Functoriality of $M^c$ on quotient stacks}
We have defined $LM^c$ as a functor from $\mathcal{A}$ to $DM(k)$. We now seek to show that it in fact induces a natural functor from the category of quotient stacks to $DM(k)$.

Consider the category $\mathcal{R}$ of resolutions of quotient stacks by algebraic spaces. 

To be precise about the definition of $\mathcal{R}$, it helps to extend the category $\mathcal{A}$ to a category $\bar{\mathcal{A}}$, which is defined just like $\mathcal{A}$, except with quotient stacks in place of algebraic spaces. We define proper and flat morphisms in $\bar{\mathcal{A}}$ just as in $\mathcal{A}$.

There is a natural embedding $\mathcal{A}\to \bar{\mathcal{A}}$. Let $\mathcal{Q}$ be the category of all quotient stacks. Then we embed $\mathcal{Q}\to \bar{\mathcal{A}}$ by constant zig-zags. 

We define $\mathcal{R}$ to be the category of factorizations $U\to V\to X$ in $\bar{\mathcal{A}}$, where $X\in \mathcal{Q}\subset \bar{\mathcal{A}}$, $U\in \mathcal{A}\subset \bar{\mathcal{A}}$, the map $V\to X$ is a termwise vector bundle, and $U\to V$ is an open embedding whose closed complement has codimension approaching infinity. This is a rephrasing of the concept of a resolution of $X$ by algebraic spaces.

We define a flat map in $\mathcal{R}$ from $U\to V\to X$ to $U'\to V'\to X'$ to be a commutative diagram
\[
\xymatrix{
	U \ar[r] \ar[d] & V \ar[r] \ar[d] & X \ar[d]\\
	U' \ar[r] & V' \ar[r] & X'
}
\]
where all the vertical maps are flat in $\mathcal{A}$. A proper map in $\mathcal{R}$ is defined by a similar commutative diagram, where the vertical maps are proper and the squares are termwise fiber squares.

Note that, as explained in the beginning of the section, every quotient stack $X\in \mathcal{Q}$ admits a resolution $U\to V\to X$ in $\mathcal{R}$.

Define the roof
\[
\xymatrix{
	& \mathcal{R} \ar[dl]_{\pi_{\mathcal{A}}} \ar[dr]^{\pi_{\mathcal{Q}}} &\\
	\mathcal{A} & & \mathcal{Q}
}
\]
by $\pi_{\mathcal{A}}(A)=U$ and $\pi_{\mathcal{Q}}(A)=X$, where $A\in \mathcal{R}$ is a factorization $U\to V\to X$ in $\mathcal{R}$.

The goal of this section is to show that the functor $LM^c\circ \pi_{\mathcal{A}}:\mathcal{R}\to DM(k)$ descends via $\pi_{\mathcal{Q}}$ to induce a functor $M^c:\mathcal{Q}\to DM(k)$. We now make this idea precise.

Define $LM^c:\mathcal{R}\to DM(k)$ by $LM^c(A)=LM^c(U)(-d_\chi)[-2d_\chi]$, where $A\in \mathcal{R}$ is a factorization $U\to V\to X$ and $d_\chi$ is the rank of $V_1$ over $X_1$. Then $LM^c$ is contravariantly functorial in flat maps in $\mathcal{R}$, and covariantly in proper maps. Similarly, $LM^c$ is functorial on fiber squares in $\mathcal{R}$ which mix proper and flat maps.

Given an object $X\in \bar{\mathcal{A}}$, let $\mathcal{R}_X$ denote the fiber of $X$ under the functor $\pi_{\mathcal{Q}}:\mathcal{R}\to \bar{\mathcal{A}}$.

\begin{construction}
	\label{canonicalisos}
	Let $X\in \mathcal{Q}$. Then given two different resolutions $A$ and $B$ in $\mathcal{R}_X$, we construct a canonical isomorphism $\varphi_{AB}:LM^c(A)\cong LM^c(B)$, such that
	
	a) We have $\varphi_{AA}=id$ for any $A\in \mathcal{R}_X$.
	
	b) We have $\varphi_{BC}\circ \varphi_{AB}=\varphi_{AC}$ for any three objects $A,B,C\in \mathcal{R}_X$.
\end{construction}
\begin{proof}
	Let $A$ be the factorization $U\to V\to X$ and let $B$ be the factorization $U'\to V'\to X$. Consider the product $A\times_{\mathcal{R}_X} B$, which is the factorization $U\times_X U'\to V\times_X V'\to X$, these fiber products being taken in $\bar{\mathcal{A}}$.
	
	Let $\pi^{AB}_A:A\times_{\mathcal{R}_X} B\to A$ and $\pi^{AB}_B:A\times_{\mathcal{R}_X} B$ be the natural projections in $\mathcal{R}_X$. 
	
	Note the projection $U\times_X U'\to U$ decomposes as $U\times_X U'\xrightarrow{i} U\times_X V'\xrightarrow{p} U$. The map $i$ is an open embedding in $\mathcal{A}$ satisfying the hypothesis of Corollary \ref{vanishingcodim}, since its closed complement is $U\times_X (V'-U')\subset U\times_X V'$, and we know that the codimension of $V'-U'$ in $V'$ approaches infinity by definition of $\mathcal{R}$. Thus $LM^c(i)$ is an isomorphism. The map $p$ consists termwise of vector bundles, so $LM^c(p)$ is an isomorphism as well. Thus $LM^c(\pi^{AB}_A)$ is an isomorphism. Symmetrically, so is $LM^c(\pi^{AB}_B)$.
	
	We then define $\varphi_{AB}=(LM^c(\pi^{AB}_B))^{-1}\circ LM^c(\pi^{AB}_A)$.
	
	We now prove the functorial conditions a),b):
	
	a) Automatic from the definition.
	
	b) Let $A,B,C\in \mathcal{R}_X$, with $A,B$ be as above and $C$ being the factorization $U''\to V''\to X$. 
	
	Let $\pi^{ABC}_A, \pi^{ABC}_B, \pi^{ABC}_C$ be the natural projections from $A\times_{\mathcal{R}_X} B \times_{\mathcal{R}_X} C$ to $A,B,C$, respectively. We proved above that a projection from a product to a factor induces an isomorphism upon applying $LM^c$, so this holds for all these projection maps. We can rewrite the formulas for $\varphi_{AB},\varphi_{BC},\varphi_{AC}$ by pulling back to $A\times_{\mathcal{R}_X} B \times_{\mathcal{R}_X} C$:
	
	We have \[\varphi_{AB}=(LM^c(\pi^{AB}_B))^{-1}\circ LM^c(\pi^{AB}_A)\]
	\[
	=(LM^c(\pi^{ABC}_B))^{-1}\circ LM^c(\pi^{ABC}_A)\]
	and similarly
	\[\varphi_{BC}= (LM^c(\pi^{ABC}_C))^{-1}\circ LM^c(\pi^{ABC}_B)\]
	and
	\[\varphi_{AC}= (LM^c(\pi^{ABC}_C))^{-1}\circ LM^c(\pi^{ABC}_A)\]
	combining the above three equations we see that indeed $\varphi_{AC}=\varphi_{BC}\circ \varphi_{AB}$.
\end{proof}
\noindent The above construction shows that, given $X\in \mathcal{Q}$, we may define $M^c(X)$ up to canonical isomorphism by the formula $M^c(X)=LM^c(A)$ for any resolution $A\in \mathcal{R}_X$.
\\
\\We record the following result, which will be used in the next section and has a similar proof.
\begin{lemma}
	\label{generalresolution}
	Consider a quotient stack $X$ and an object $A\in \mathcal{A}$ together with a flat morphism $A\to X$ in $\mathcal{A}$ of relative dimension $r$. Suppose that the flat pullback $M^c(Y)\to LM^c(A\times_X Y)(-r)[-2r]$ is an isomorphism for all algebraic spaces $Y$ over $X$. Then so is the flat pullback $M^c(X)\to LM^c(A)$. 
\end{lemma}
\begin{proof}
	For any representable morphism of stacks $Y\to X$, let $H(Y)\in Sp_{Tate}(k)$ denote $LM^c(A\times_X Y)$. Let $\bar{\mathcal{A}}_X$ be the category of objects in $\bar{\mathcal{A}}$ equipped with a representable map to $X$. Then  $H$ can be extended to a functor $\hat{H}:\bar{\mathcal{A}}_X\to Sp_{Tate}(k)^{ZigZag}$ by evaluating $H$ at each morphism of $\bar{\mathcal{A}}$. We can then define $LH:\bar{\mathcal{A}}_X\to DM(k)$ by $holim\circ \hat{H}$. 
	
	Now, fix a resolution $U\to V\to X$ of $X$ by algebraic spaces. Consider the sequence
	\[
	LH(V-U)\to LH(V)\to LH(U)
	\]
	Note that $H$ satisfies a localization triangle $H(Y-U)\to H(Y)\to H(U)$ for any open embedding $U\to Y$ of stacks, by Lemma \ref{holimtriangles}. Similarly, the functor $LH$ satisfies a localization triangle $LH(A-U)\to LH(A)\to LH(U)$, for any open embedding $U\to A$ in $\bar{\mathcal{A}}_X$, by another application of Lemma \ref{holimtriangles}. Thus the above sequence is an exact triangle in $DM(k)$.
	
	We know $U$ consists of algebraic spaces, and by assumption $H(Y)$ is naturally isomorphic to $M^c(Y)$ for any algebraic space $Y$ over $X$. Thus 
	\[LH(U)\cong LM^c(U)(-r)[-2r]\cong M^c(X)\]
	by definition of the compactly supported motive of a quotient stack.
	
	On the other hand, the vector bundle $V\to X$ induces an isomorphism of motives $LH(V)\cong H(X)=LM^c(A)$. 
	
	Finally, we claim $LH(V-U)=0$. Indeed, for each term $V_i-U_i$ in $V-U$, we see that $H(V_i-U_i)$ is a homotopy limit of motives in $D_{-r_i}$, where $r_i$ is the codimension of $V_i-U_i$ in $V_)$. Thus $H(V_i-U_i)$ is in $E_{-r_i}^{\bot}$, since $D_{-r_i}\subset E_{-r_i}^{\bot}$ and $E_{-r_i}^{\bot}$ is colocalizing. 
	
	Therefore $LH(V-U)=holim H(V_i-U_i)$ is in the intersection of the $E_{-r_j}^{\bot}$, which is $0$ since $r_j\to \infty$.
	
	In conclusion, since the homotopy limits of the three sequences form an exact triangle, we see that the map $LM^c(A)\to M^c(X)$ is an isomorphism, as desired.
\end{proof}
Having shown that the motive $M^c(X)$ is defined up to canonical isomorphism, we now turn to its functoriality, studying descent of morphisms under the functor $\pi_{\mathcal{Q}}:\mathcal{R}\to \mathcal{Q}$.
\begin{proposition}
	\label{morphismdescent}
	Let $f_1:A\to B$ and $f:A'\to B'$ be (flat or proper) morphisms in $\mathcal{R}$ such that $\pi_{\mathcal{Q}}(f_1)=\pi_{\mathcal{Q}}(f_2)$. Then $LM^c(f_1)$ and $LM^c(f_2)$ are conjugate via the canonical isomorphisms $\varphi_{AA'}$ and $\varphi_{BB'}$.
\end{proposition}
\begin{proof}
	Assume first that $f_1,f_2$ are flat maps. Let $f:A\times_{\mathcal{R}_X} A'\to B\times_{\mathcal{R}_Y} B'$ be the natural map. Consider the following square in $\mathcal{R}$:
	\[
	\xymatrix{
		A\times_{\mathcal{R}_X} A' \ar[r]^{f} \ar[d]^{\pi_A} & B\times_{\mathcal{R}_Y} B' \ar[d]^{\pi_B}\\
		A \ar[r]^{f_1} & B
	}
	\]
	By functoriality of $LM^c$ on $\mathcal{R}$ with respect to flat maps, we have
	\[
	LM^c(f)=LM^c(\pi_A)\circ LM^c(f_1)\circ LM^c(\pi_B)^{-1}
	\]
	By similarly considering projections to $A'$ and $B'$, we have
	\[
	LM^c(f)=LM^c(\pi_{A'})\circ LM^c(f_2)\circ LM^c(\pi_{B'})^{-1}
	\]
	Combining the two equations, we see
	\[
	LM^c(\pi_{A'})^{-1}\circ LM^c(\pi_A)\circ LM^c(f_1)=LM^c(f_2)\circ LM^c(\pi_{B'})^{-1}\circ LM^c(\pi_B)
	\]
	which by the details of Construction \ref{canonicalisos} amounts to
	\[
	\varphi_{AA'}\circ LM^c(f_1)=LM^c(f_2)\circ \varphi_{BB'}
	\]
	as desired.
	\\
	\\Now assume that $f_1,f_2$ are proper maps. Note that the above square becomes a termwise fiber square, so by functoriality of $LM^c$ on $\mathcal{R}$ with respect to fiber squares of flat and proper maps, we see by a similar argument that
	\[
	\varphi_{BB'}\circ LM^c(f_1)=LM^c(f_2)\circ \varphi_{AA'}
	\]as desired.
\end{proof}
\noindent Note that any flat map $f:X\to Y$ in $\mathcal{Q}$ is covered by a flat map $\bar{f}:A\to B$ in $\mathcal{R}$, in the sense that $\pi_{\mathcal{Q}}(\bar{f})=f$. Indeed, first fix a resolution $A_0\in \mathcal{R}_X$ and $B\in \mathcal{R}_Y$. Then define $A=f^*B\times_X A_0$. That is, writing $A_0$ as the composition $U\to V\to X$ and writing $B$ as the composition $O\to W\to Y$, we define $A\in \mathcal{R}_X$ as the composition $U\times_X f^*O\to V\times_X f^*W\to X$. In general, $f^*O$ might not be in $\mathcal{A}$ (that is, a zig-zag of algebraic spaces), but we know $U\times_X f^*O\in \mathcal{A}$, since it is an open embedding of a termwise vector bundle over $U\in \mathcal{A}$. The composition $A\to f^*B\to B$ covers the map $f$.

On the other hand, not every proper map $f:X\to Y$ in $\mathcal{Q}$ can be covered by a proper map $\bar{f}:A\to B$ in $\mathcal{R}$. The point is that, by our definition of proper maps in $\mathcal{R}$, the object $A$ would have to be the termwise pullback $f^*B$, which may not always be in $\mathcal{R}$. For example, if $G$ is a finite group, then $BG\to Spec(k)$ is proper, but given any resolution $U\to V\to Spec(k)$ in $\mathcal{R}$, the pullback $f^*U$ equals the product $BG\times_k U$, which is not a zig-zag of algebraic spaces. 

For this reason, we restrict to proper maps in $\mathcal{Q}$ which are termwise representable. Such maps, which in the future we simply refer to as representable, can always be covered by proper maps in $\mathcal{R}$.
\begin{theorem}
	\label{mcfunctor}
	There is a functor $M^c\colon \mathcal{Q}\to DM(k)$, canonical up to natural isomorphism, which is contravariantly functorial with respect to flat maps and covariantly functorial with respect to representable proper maps, and functorial with respect to fiber squares mixing flat maps and representable proper maps.
\end{theorem}
\begin{proof}
	For each object $X\in \mathcal{Q}$, fix a resolution $A_X\in \mathcal{R}_X$. Then define 
	\[M^c(X):=LM^c(A_X)\]
	Similarly, for each morphism $f:X\to Y$ in $\mathcal{Q}$ which is either flat or representable proper, fix a cover $\bar{f}:A_{f,X}\to A_{Y}$ of $f$ in $\mathcal{R}$. Then define the morphism $M^c(f):LM^c(Y)\to LM^c(X)$ by
	\[
	M^c(f)=\varphi_{A_{f,X}A_X}\circ LM^c(\bar{f}) \text{ if }f \text{ is a flat map}
	\]
	\[
	M^c(f)=LM^c(\bar{f})\circ \varphi_{A_XA_{f,X}} \text{ if }f \text{ is a representable proper map}
	\]
	
	To show functoriality, consider a sequence of flat maps $X\xrightarrow{f}Y\xrightarrow{g}Z$ in $\mathcal{Q}$. The key is that there is a cover $A\to B\to C$ of in $\mathcal{R}$ of $X\to Y\to Z$ consisting of flat maps. It then follows formally from functoriality of $LM^c$ on $\mathcal{R}$, together with the results of Construction \ref{canonicalisos} and Proposition \ref{morphismdescent}, that $M^c(f)\circ M^c(g)=M^c(g\circ f)$.
	
	Functoriality with respect to representable proper maps is proven in the same way.
	
	Functoriality with respect to fiber squares mixing flat and representable proper maps is also proven similarly, by noting that given any fiber square
	\[
	\xymatrix{
		X'\ar[r] \ar[d] & Y' \ar[d]\\
		X\ar[r] & Y
	}
	\]
	in $\mathcal{Q}$ with vertical maps flat and horizontal maps representable proper maps, this square is covered by some fiber square in $\mathcal{R}$ similarly of the form
	\[
	\xymatrix{
		A'\ar[r] \ar[d] & B' \ar[d]\\
		A\ar[r] & B
	}
	\]
	where the vertical maps are flat and the horizontal maps are proper.
	
	Finally, we note that a different choice of representatives $A_X$ and $A_{f,X}$ would yield a naturally isomorphic functor via the isomorphisms in Construction \ref{canonicalisos}.
\end{proof}
\begin{corollary}
	\label{groupfunctoriality}
	The association $GpSch/k\to DM(k)$ given by $G\mapsto M^c(BG)(\text{dim}G)[2\text{dim}G]$ is contravariantly functorial.
\end{corollary}
\begin{proof}
	Given a group homomorphism $H\to G$, there is a flat morphism of quotient stacks $BH\to BG$.
\end{proof}
\chapter{The Mixed Tate Property}
Throughout this section, we assume that we have fixed a field $k$ and a coefficient ring $R$ such that the exponential characteristic of $k$ is invertible in $R$. All statements about motives are understood to be with coefficients in $R$.
\section{Background}
We describe three different types of Kunneth properties, from weaker to stronger.
\begin{definition}
	\label{kunneth}
	Let $X$ be a quotient stack of finite type over a field $k$.
	
	$X$ has the \textit{weak Kunneth property} if for any field extension $K$ of $k$, the natural morphism $CH_*(X)\to CH_*(X_K)$ is a surjection.
	
	$X$ has the \textit{Chow Kunneth property} if for any separated scheme $Y$ of finite type over $k$, the natural morphism
	\[
	CH_*(X)\otimes CH_*(Y)\to CH_*(X\times_k Y)
	\]
	is an isomorphism.
	
	$X$ has the \textit{motivic Kunneth property} if for any $k$-scheme $Y$, the Kunneth spectral sequence for $X$ and $Y$ converges to motivic homology groups of $X\times Y$.\cite{burt} 
\end{definition}
The motivic Kunneth property implies the Chow Kunneth property, which in turn implies the weak Chow Kunneth property. Totaro characterized these properties in terms of the motive $M^c(X)\in DM(k;R)$.

The triangulated category $DMT(k;R)$ of \textit{mixed Tate motives} is the smallest localizing subcategory of $DM(k;R)$ that contains the Tate motives $R(j)$ for all integers $j$.\cite{rondigs}
\begin{proposition} (Totaro\cite{burt})
	The quotient stack $X$ satisfies the motivic Kunneth property if and only if the motive $M^c(X)$ is a mixed Tate motive.
\end{proposition}
Note that, since $DMT(k;R)$ is triangulated, if we are given an exact triangle $A\to B\to C$ in $DM(k;R)$ such that two out of three of the motives are in $DMT(k;R)$, then so is the third. For example, this may be applied to the localization triangle of $M^c(Z)\to M^c(X)\to M^c(X-Z)$ for any substack $Z$ of a quotient stack $X$, showing that if two of $Z,X,X-Z$ have mixed Tate motives, so does the third.

A related notion is a linear scheme. The family of \textit{linear schemes} over $k$ is the smallest set of schemes containing affine space $\mathbb{A}^n_k$ for all $n\ge 0$ and such that for any scheme $X$ of finite type over $k$ with a closed subscheme $Z$, if $Z$ and $X-Z$ are linear schemes, then so is $X$, and if $X$ and $Z$ are linear schemes, then so is $X-Z$. By the definition of mixed Tate motives, together with the localization exact triangle, it follows by induction on dimension that whenever $X$ is a linear scheme, the motive $M^c(X)$
is mixed Tate.

We say that an affine group scheme $G$ satisfies the \textit{mixed Tate property} if $M^c(BG)$ is mixed Tate, which, by the above Proposition, is equivalent to $BG$ satisfying the motivic Kunneth property.

The following result makes the mixed Tate property easier to check.
\begin{proposition} (Totaro\cite{burt})
	Let $X$ be a quotient stack of finite type over $k$ and let $E\to X$ be a $GL(n)$-bundle. Then $M^c(X)$ is mixed Tate if and only if $M^c(E)$ is mixed Tate.
\end{proposition}
\begin{corollary}
	Let $G\subset GL(n)$ be a subgroup scheme. Then $G$ has the mixed Tate property if and only if $M^c(GL(n)/G)$ is mixed Tate.
\end{corollary}
We note basic results about the mixed Tate property.
\begin{lemma}
	\label{mixedtatefunctoriality}
	Let $G,H$ denote affine group schemes.
	\\
	a) The product $G\times H$ has the mixed Tate property if and only if both $G$ and $H$ have the mixed Tate property.
	\\b)(Totaro)\cite{burt} If $BG$ can be approximated by linear schemes, then so can the wreath product $\mathbb{Z}/p\wr G$.
	\\b') In particular, if $G\subset GL(n)$ and $GL(n)/G$ is a linear scheme, then $BG$ can be approximated by linear schemes, so $\mathbb{Z}/p\wr G$ has the mixed Tate property.
\end{lemma}
\begin{proof}
	a) Suppose $G$ and $H$ have the mixed Tate property. Then $B(G\times H)\cong BG\times BH$ does as well, sinde $DM(k;R)$ is closed under tensor product.
	
	Conversely, suppose $G\times H$ has the mixed Tate property. We may form the sequence $G\to G\times H\to G$, which induces by functoriality a sequence 
	\[M^c(BG)\to M^c(B(G\times H))(\text{dim}H)[2\text{dim}H] \to M^c(BG) \]
	whose composition is the identity, demonstrating that $M^c(BG)$ is a summand of $M^c(B(G\times H))$, so it has the mixed Tate property as well since the category of mixed Tate motives is localizing.
	
	b) See \cite{burt}.
	
	b') Let $U_N$ consist of linearly independent $n$-tuples in $\mathbb{A}^{N+n}$. Then the projection $\pi:U_N\to Gr(n,n+N)$ mapping a tuple to its span is a $GL(n)$-bundle. The variety $Gr(n,n+N)$ has an affine stratification into cells such that $\pi$ is a trivial $GL(n)$-bundle over each cell. Thus the quotient $p:U_N/G\to U_N/GL(n)\cong Gr(n,n+N)$ is a $GL(n)/G$-bundle which is trivial over each cell. Thus $U_N/G$ can be stratified into pieces of the form $\mathbb{A}^m\times GL(n)/G$, which are linear if $GL(n)/G$ is linear, showing that $U_N/G$ is linear, so $BG$ may be approximated by linear schemes.
\end{proof}
\section{Techniques}
\begin{lemma}
	\label{orbits}
	Let $G$ be a smooth affine algebraic group over $k$ acting on a scheme $X$ over $k$. Suppose that the induced action of $G(\bar{k})$ on $X(\bar{k})$ has finitely many orbits $O_i\subset X(\bar{k})$ and that, furthermore, each orbit $O_i$ contains a $k$-point  $x_i$ (that is, a point of $X(\bar{k})$ which comes from $X(k)$). Then $X$ can be stratified by quotient schemes $G/H_i$, where $H_i\subset G$ is the stabilizer group of $x_i$. That is, there is a chain of reduced closed subschemes $X=X_0\supset X_1 \supset \cdots \supset X_n=\emptyset$ such that, for an appropriate reordering of the $x_i$, $G/H_i$ is isomorphic to an open subscheme $U_i\subset X_{i-1}$, whose complement is $X_i$.
\end{lemma}
\begin{proof}
	Since the action of $G(\bar{k})$ on $X(\bar{k})$ has finitely many orbits, each irreducible component of $X(\bar{k})$ must be contained in (exactly) one of the closures $\overline{O_i}$. So fix an irreducible component $V\subset X$ and an orbit $O_1$ such that $\overline{O_1}$ contains $V$. Note that $O_1$ does not intersect any $\overline{O_i}$ for $i\ne 1$, since if it did, then $\overline{O_i}$, being stable under $G$, must contain all of $O_1$, a contradiction. Thus $O_1$ must be contained in the open complement of $\bigcup_{i\ne 1} \overline{O_i}$, so it equals this open set, demonstrating that $O_1$ is open in $X(\bar{k})$.
	
	Now, by our assumption that each orbit contains a $k$-point, let $x_1\in O_1$ be a $k$-point. Let $O_{x_1}\subset X$ be the orbit scheme of $x_1$. So we know that $O_{x_1}(\bar{k})=O_1$ and the orbit $O_{x_1}$ is open in $X$. Also, since $x_1$ is a $k$-point, the orbit map $G\to O_{x_1}: g\mapsto gx$ gives $O_{x_1}$ the structure of a quotient $G/H_1$, where $H_1$ is the stabilizer group of $x_1$. Now, the closed complement $X\backslash O_{x_1}$, with the reduced scheme structure, is stable under $G$. The orbits of the action of $G(\bar{k})$ on $X(\bar{k})\backslash O_1$ are exactly the orbits of $X(\bar{k})$ other than $O_1$, and the stabilizer of a $k$-point  $x\in X\backslash O_{x_1}$ is isomorphic to the stabilizer of $x\in X$, so we may now apply induction.
\end{proof}
The following corollaries are two common situations:
\begin{corollary} 
	\label{GLnaction}
	Let $X$ be a mixed Tate scheme with a $GL(n)$ action. Suppose that $GL(n)(\bar{k})$ acts on $X(\bar{k})$ with finitely many orbits each of which contains a $k$-point. Let $x_i$ be $k$-points which are representatives of all the orbits of $GL(n)(\bar{k})$ on $X(\bar{k})$. For each $i$, let $H_i\subset GL(n)$ denote the stabilizer group of $x_i$. Assume the orbit of $x_1$ is an open orbit. Suppose that the groups $H_i$ satisfy the mixed Tate property for all $i>1$. Then $H_1$ also satisfies the mixed Tate property.
\end{corollary}
\begin{proof}
	Assume without loss of generality that the $x_i$ are ordered as in the statement of lemma 1.1 and let $X_i,U_i$ be as in the statement of lemma 1.1. We prove by induction that each $X_i$, $i\ge 1$, is mixed Tate (the base case being $X_n=\emptyset$). Indeed, let $i\ge 1$ and assume $X_{i+1}$ is mixed Tate. By the assumption that each $H_{i+1}$ has the mixed Tate property ,and since $U_{i+1}\cong GL(n)/H_{i+1}$, the open subscheme $U_{i+1}\subset X_{i}$ is mixed Tate. But we assumed $X_{i+1}=X_{i}\backslash U_{i+1}$ is mixed Tate, so it follows that $X_{i}$ is mixed Tate, as desired. Since $X_1=X\backslash U_1$ is mixed Tate and $X$ is mixed Tate, it follows that the open subscheme $U_1\subset X$ is mixed Tate. Since $U_1\cong GL(n)/H_1$, the group $H_1$ has the mixed Tate property.
\end{proof}
\begin{corollary}
	\label{solvableaction}
	Let $G$ be a split solvable algebraic group over $k$ acting on a scheme $X$. Suppose $G(\bar{k})$ acts on $X(\bar{k})$ with finitely many orbits, each of which has a $k$-point. Then $X$ is mixed Tate.
\end{corollary}
\begin{proof}
	Let $X_i, U_i$ be as in the statement of lemma 1.1. By induction, we assume $X_1$ is mixed Tate. Also $U_1\cong G/H_1$ is mixed Tate since all homogeneous spaces of split solvable groups are mixed Tate. Since $X_1$ is the complement of the open subscheme $U_1\subset X$, we conclude that $X$ is mixed Tate.
\end{proof}
Note: the assumption that every orbit has a $k$-point is necessary for the above corollary. As a simple example, let $G=1$ act on $X=\text{Spec} \mathbb{C}$ over $\mathbb{R}$. This action has two orbits over $\mathbb{C}$ since $\text{Spec}\mathbb{C}\times_{\mathbb{R}} \text{Spec}\mathbb{C}$ consists of two points. But $X=\mathbb{C}$ is not mixed Tate over $\mathbb{R}$ since it does not have the weak Chow Kunneth property.
\noindent We now describe a common type of group which has the mixed Tate property. Given a sequence $n_1,\dots, n_k$, set $n=n_1+\cdots +n_k$. Let $T_{n_1,\dots,n_k}\subset GL(n)$ denote the stabilizer of the flag $F_1\subset \cdots \subset F_k=V$, where $V$ is an $n$-dimensional vector space with basis $e_i$ and $F_i=<e_1,\dots, e_{n_1+\cdots+n_i}>$. Concretely, $T_{n_1,\dots,n_k}$ consists of block-upper-triangular matrices with blocks of sizes $n_1,\dots, n_k$ along the diagonal. Let $T'_{n_1,\dots, n_k}$ be the transpose of $T_{n_1,\dots, n_k}$, that is, consisting of block-lower-triangular matrices of the same sizes.

Let $D_{n_1,\dots, n_k}\cong GL(n_1)\times \cdots \times GL(n_k)\subset GL(n)$ denote the subgroup of $GL(n)$ stabilizing all of the subspaces $<e_1,\dots, e_{n_1}>,\dots, <e_{n_{k-1}+1},\dots, e_{n_k}>$. Concretely, $D_{n_1,\dots, n_k}$ consists of block-diagonal matrices with blocks of sizes $n_i$. We have a group homomorphism $\pi: T_{n_1,\dots, n_k}\to D_{n_1,\dots, n_k}$ given by picking out the diagonal terms. Let $U_{n_1,\dots, n_k}$ denote the kernel of $\pi$, consisting of block-unipotent matrices. Then $U_{n_1,\dots, n_k}$ is isomorphic as a variety to affine space $\mathbb{A}^N$, where $N=\sum_{i<j} n_in_j$. Similarly, we have a map $\pi':T_{n_1,\dots, n_k}\to D_{n_1,\dots, n_k}$ and we denote its kernel $U'_{n_1,\dots,n_k}$.
\begin{lemma}
	\label{commonexample}
	We use the above notation. Suppose $H\subset GL(n)$ can be expressed as the semi-direct product $J\rtimes K$, where $K\subset D_{n_1,\dots, n_k}$ and $J\subset U_{n_1,\dots, n_k}$ (respectively, $J\subset U'_{n_1,\dots, n_k}$), where the multiplication operation on $J\rtimes K$ comes from $GL(n)$. Suppose $K$ has the mixed Tate property and $J$ is a linear subspace of $U_{n_1,\dots, n_k}\cong \mathbb{A}^N$ (resp. of $U'_{n_1,\dots, n_k}$). Then $H$ has the mixed Tate property.
\end{lemma}
\begin{proof} By assumption, $H=J\rtimes K$, so every element $h\in H$ can be written (uniquely) as $h=j\cdot k$, with $k\in K$ and $j\in J$. In fact, we must have $k=\pi(h), j=hk^{-1}$. Thus the map $\varphi:H\to J$ given by $h\mapsto h\pi(h)^{-1}$  induces an isomorphism of varieties $\bar{\varphi}:H/K\to J$. The group $H$ acts on $H/K\cong J\cong \mathbb{A}^m$ via $\bar{\varphi}$. We verify that $H$ acts on $\mathbb{A}^m$ by affine-linear transformations. It suffices to verify that both $K$ and $J$ act by affine linear transformations. Indeed,   Indeed, consider elements $h\in H$ and $j\in J$. Write $h=kj;$, where $k\in K$, $j'\in J$. Then $\bar{h}\cdot u = j'u\in J$. Note that multiplication by a fixed element of $U_{n_1,\dots, n_k}$ is affine-linear on $U_{n_1,\dots, n_k}\cong \mathbb{A}^N$. Since $J\subset U_{n_1,\dots, n_k}$ is a linear subspace, the induced map on $J$ is also affine-linear. Thus $H$ acts on $J\cong \mathbb{A}^m$ by affine-linear transformations, as desired.
	
	An algebraic group $G$ is \textit{special} if every $G$-bundle is Zariski-locally trivial. For example, the group of affine transformations is \textit{special}.
	
	Now, there is a natural map $p: GL(n)/K\to GL(n)/H$ with fiber $H/K\cong J\cong \mathbb{A}^{m}$. Since $H$ acts by affine transformations on $J$ and the group of affine transformations is special, $p$ is Zariski-locally trivial. Since $K$ has the mixed Tate property, $GL(n)/K$ is mixed Tate, hence, since $p$ is an affine-space bundle, $GL(n)/H$ is mixed Tate, so $H$ has the mixed Tate property.
	
	The case of $J\subset U'_{n_1,\dots, n_k}$ instead of $U_{n_1,\dots, n_k}$ is similar.
\end{proof}
\section{$Sp(2n), O(q)$ and $SO(q)$}
Before treating the mixed Tate property for these groups, we first record the known computations of the Chow groups of their classifying spaces. Note these Chow group calculations hold for any field, already indicating the weak Chow Kunneth property.
\begin{proposition}(Totaro\cite{burtchow})
	We have the following computations
	\label{chowclassical}
	\[
	CH^*BGL(n)=\mathbb{Z}[c_1,\dots,c_n]
	\]
	\[
	CH^*BO(n)=\mathbb{Z}[c_1,\dots,c_n]/(2c_i=0 \text{ for i odd})
	\]
	\[
	CH^*BSp(2n)=\mathbb{Z}[c_2,c_4,\dots, c_{2n}]
	\]
	\[
	CH^*BSO(2n+1)=\mathbb{Z}[c_2,c_3\dots,c_{2n+1}]/(2c_i=0 \text{ for i odd})
	\]
\end{proposition}
Here $c_i$ stand for the chern classes under the natural representations.

The Chow groups of $SO(2n)$ involve a more difficult calculation. Indeed they stand out for not being generated by Chern classes.
\begin{proposition}(Field\cite{field})
	\[
	CH^*BSO(2n)= \mathbb{Z}[c_2,c_3\dots,c_{2n},y_n]/(2c_{odd},y_n\cdot c_{odd},y_n^2+(-1)^n2^{2n-2}c_{2n})
	\]
\end{proposition}
\begin{proposition}
	\label{Sp2n}
	Over an arbitrary field $k$ and for any $n>0$, the symplectic group $Sp(2n)$ over $k$ has the mixed Tate property.
\end{proposition}
\begin{proof}
	We consider the action of $GL(2n)$ on the space $V$ of alternating forms of dimension $2n$. Let $f_0=0$ and let $f_i=x_1\wedge x_2 + \cdots x_{2i-1}\wedge x_{2i)}$, for $i=1,\dots,n$, and $f_{n+1}=0$. Let $H_i\subset GL(2n)$ be the stabilizer of $f_i$. Then $H_i=K_i\rtimes J_i$, where $K_i\subset D_{2i,2n-2i}$ consists of matrices of the block-matrix form $\left[\begin{array}{cc}
	A & 0
	\\0 & B
	\end{array}\right]$ and $J_i= U'_{2i,2n-2i}$. By induction, if $i<n$ then $Sp(2i)$ has the mixed property, so $K_i$ has the mixed Tate property. So by lemma \ref{commonexample}, $H_i$ has the mixed Tate property. The forms $f_i$ form a set of representatives for all orbits of $GL(n)(k)$ on $V(k)$, and the orbit of $f_n$ is the open orbit of nondegenerate alternating forms. Since all the stabilizers $H_i$ of the $f_i$ with $i<n$ have the mixed Tate property, it follows by Corollary \ref{GLnaction} that $Sp(2n)$, the stabilizer of $f_n$, has the mixed Tate property, as desired.
\end{proof}
For the remainder of this section, let $k$ be a field of characteristic not equal to $2$. We now show that $O(q)$ has the mixed Tate property for any nondegenerate quadratic form $q$ of dimension $n$. We could use the same idea as above, analyzing the $GL(n)$-orbits on the space of quadratic forms, but instead we take a different approach, analyzing the $B$-orbits on the homogeneous space $GL(n)/O(q)$, where $B\subset GL(n)$ is the Borel subgroup consisting of upper-triangular matrices. This approach lends itself better to the study of $SO(q)$. We identify the homogeneous space $GL(n)/O(q)$ with the space $P_n$ of nondegenerate quadratic forms $q$ of dimension $n$. The action of $GL(n)$ is given by $(A\cdot q)(v)=q(Av)$ and the quadratic form $q\in P_n$ is taken as the neutral element of $GL(n)/O(q)$. We also identify $P_n$ with the space of nonsingular symmetric matrices, by associating to each quadratic form $q'$ the unique symmetric matrix $Q'\in GL(n)$ such that $q'(v)=v^TQ'v$.
\begin{lemma}
	\label{GLnOnBorbits}
	Let $S_n$ be the symmetric group on $n$ letters and let $I_n\subset S_n$ consist of all involutions. Then the permutation matrices corresponing to elements of $I_n$ form a set of representatives for the orbits of $B(\bar{k})$ on $(GL(n)/O(n))(\bar{k})$.
\end{lemma}
\begin{proof}
	See \cite{richspring}.
\end{proof}
\begin{proposition}
	\label{Oq}
	Let $k$ be an arbitrary field of characteristic not equal to $2$ and let $q$ be a nondegenerate quadratic form over $k$. Then the orthogonal group $O(q)$ has the mixed Tate property.
\end{proposition}
\begin{proof}
	The above lemma shows that $B(\bar{k})$ acts on $(GL(n)/O(n))(\bar{k})$ with finitely many orbits, each of which has a point defined over $k$ (since permutation matrices are defined over $\mathbb{Z}$). By corollary $1.3$, it follows that $GL(n)/O(q)$ has the mixed Tate property.
\end{proof}
We now study the mixed Tate property for $SO(q)$, for a nondegenerate quadratic form $q$. If $q$ has dimension $2n+1$, this is easy.
\begin{proposition}
	\label{SOqodd}
	Let $k$ be an arbitrary field of characteristic not equal to $2$ and let $q$ be a nondegenerate quadratic form over $k$ of dimension $2n+1$. Then the group $SO(q)$ has the mixed Tate property.
\end{proposition}
\begin{proof}
	We have $O(q)\cong SO(q)\times \mu_2$. Since $O(q)$ has the mixed Tate property, so does $SO(q)$ by $1.4a)$.
\end{proof}
The even-dimensional case is harder. We again consider consider the action of $B$ on $GL(n)/SO(q)$. Given $n\ge 0$, let $P_{2n}\subset GL(2n)$ denote as above the space of nonsingular quadratic forms of dimension $2n$. Given $d\in k$, let $W^d_{n}$ denote the subscheme of $P_{2n}\times G_m$ consisting of pairs $(q,t)$ such that $det(q)=t^2d$. Given a quadratic form $q$ of determinant $d$, we can view $W^d_n$ as the quotient scheme $GL(2n)/SO(q)$, where $GL(2n)$ acts on $W^d_n$ by $g\cdot(q',t)=(g\cdot q', det(g)t)$ and $(q,1)$ is the neutral element of $GL(2n)/SO(q)$.
\begin{proposition}
	\label{SOqeven}
	Let $k$ be an arbitrary field of characteristic not equal to $2$ and let $q$ be a nondegenerate quadratic form over $k$ of dimension $2n$. Then:
	\\i)If $\text{det}(q) = (-1)^n (\text{mod }  (k^\times)^2)$, then the special orthogonal group $SO(q)$ has the mixed Tate property.
	\\ii)Conversely, if $\text{det}(q) \ne (-1)^n (\text{mod }  (k^\times)^2)$, then $SO(q)$ does not satisfy the weak Chow Kunneth property, so in particular it does not satisfy the mixed Tate property.
\end{proposition}
In order to prove part $ii)$ of the proposition, we  calculate $CH_*(GL(2n)/SO(q))$. Field\cite{field} calculated this over the complex numbers, though her proof works as well over any algebraically closed field of characteristic not equal to $2$. We generalize Field's computation to an arbitrary quadratic form $q$ over a general field $k$.
\begin{lemma}
	\label{SOqchow}
	i)If $\text{det}(q) = (-1)^n (\text{mod }  (k^\times)^2)$, then $CH^*(GL(2n)/SO(q))\cong \mathbb{Z}\oplus \mathbb{Z}y$, where $y$ is a codimension $n$ cycle.
	\\ii)If $\text{det}(q) \ne (-1)^n (\text{mod }  (k^\times)^2)$, then $CH^*(GL(2n)/SO(q))\cong \mathbb{Z}$.
\end{lemma}
\noindent
\textit{Proof of Proposition \ref{SOqeven} and Lemma \ref{SOqchow}}:
\\
\ref{SOqeven}i) Consider the action of the Borel subgroup $B$ of $GL(2n)$ on $GL(2n)/SO(q)$. For an element $q'\in GL(2n)/O(q)= P_{2n}$, the fiber of $q'$ under the projection $GL(2n)/SO(q)\to GL(2n)/O(n)$ consists of two elements over $\bar{k}$: $(q',\pm \sqrt{det(q')/d})$.

Let $h\in I_{2n}\subset S_{2n}\subset GL(2n)$ be the representative of $q'$ under the action of $B$ on $P_{2n}$. First suppose that $h(e_i)=e_i$ for some basis vector $e_i$. Then consider the diagonal matrix $b\in B$ with $b(e_j)=e_j$ for $j\ne i$, $b(e_i)=\sqrt{d/det(h)}e_i$. Thus $b\cdot (h,\pm \sqrt{det(h)/d})=(b^2\cdot h, \pm 1)$, which is a $k$-point. Thus $x$ has a $k$-point in its orbit. 

Now suppose that $h$ is fixed-point free as a permutation. Then $det(h)=(-1)^n$ since it is a product of $n$ transpositions. Then the condition that $d = (-1)^n \text{ (mod } (k^\times)^2)$ is precisely the condition that $\sqrt{det(h)/d}\in k$, thus implying that $x$ has a $k$-point $(h,\sqrt{det(h)/d})$ in its orbit. Thus, if $d = (-1)^n \text{mod } (k^\times)^2$, then every orbit has a $k$-point so $GL(2n)/SO(q)$ is mixed Tate by corollary \ref{solvableaction}.
\\
\\\ref{SOqchow}i) Field\cite{field} proved this result for algebraically closed fields. By \ref{SOqeven}i), if $det(q) \ne (-1)^n \text{mod } (k^\times)^2$, then $SO(q)$ has the mixed Tate property, hence by the Kunneth property, $CH^*((GL(2n)/SO(q))_k)\cong CH^*((GL(2n)/SO(q))_{\bar{k}})$.
\\
\\\ref{SOqchow}ii) Now suppose that $d:=det(q) \ne (-1)^n \text{mod } (k^\times)^2$. We closely follow the proof method of \cite{field} to show that $CH^*(GL(2n)/SO(q))=CH^*(W^d_n)\cong \mathbb{Z}$, arguing by induction. The base case is taken to be $n=0$, in which case $W^d_0\subset G_m$ is $\text{Spec}(k[t]/(t^2=\frac{1}{d}))$. In the case of $n=0$, $d \ne (-1)^n \text{mod } (k^\times)^2$ implies $\frac{1}{d}$ is not in $(k^\times)^2$. Thus, $k[t]/(t^2=\frac{1}{d})$ is a field, so $\text{Spec}(k[t]/(t^2=\frac{1}{d}))$ is irreducible, and $CH^*(W^d_0)\cong \mathbb{Z}$ as desired.  

Now we consider general $n>0$. For each $i$ between $1$ and $2n$, let $X_i\subset W^d_n$ be the subvariety consisting of pairs $(q',t)$ such that $q'(e_i,e_{2n})\ne 0$ and $q'(e_j,e_{2n})=0$ for $i< j \le 2n$. Let $\bar{q'}$ be the restriction of $q'$ to the orthogonal complement of $k<e_i,e_{2n}>$ with respect to $q'$. Then $\bar{q'}$ is a quadratic form of dimension $2(n-1)$ and we have $det(q')=det(\bar{q'})\cdot -q'(e_i,e_{2n})^2$. We can therefore define a map $X_i\to W^{-d}_{n-1}$ by $(q',t)\mapsto (\bar{q'},t/q'(e_i,e_{2n}))$. This map is a trivial fibration, with fiber $G_m\times \mathbb{A}_k^{2n+i-2}$. Since, by assumption, $d = (-1)^{n} (\text{mod }(k^\times)^2)$, it follows that $-d = (-1)^{n-1} (\text{mod }(k^\times)^2$) Thus we may apply induction to conclude $CH^*(X_i)\cong \mathbb{Z}$.

We now show that $CH^*(\overline{X_i})\cong \mathbb{Z}$ for all $i$, by induction. The base case is $i=1$, where it follows from above, since $X_1$ is closed in $W^d_n$. For each $i>1$, $X_i$ is an open subvariety of $\overline{X_i}$ whose closed complement is $\overline{X_{i-1}}$. By induction, $CH^*(\overline{X_{i-1}})\cong \mathbb{Z}$, generated by the codimension-zero cycle. Since $\overline{X_{i-1}}$ is cut out of $\overline{X_i}$ by $q'(e_i,e_{2n})=0$, we see that the pushforward $CH^*(\overline{X_{i-1}})\to CH^*(\overline{X_{i}})$ is the $0$ map, so by the basic exact sequence of Chow groups, $CH^*(\overline{X_{i}})\cong CH^*(X_i)\cong \mathbb{Z}$. In particular, since $\overline{X_{2n}}=W^d_n$, we have $CH^*(W^d_n)\cong \mathbb{Z}$, as desired. 
\\
\\\ref{SOqeven}ii) If $d:=det(q) \ne (-1)^n \text{mod } (k^\times)^2$, then it follows from Lemma \ref{SOqchow}i) and ii) that the map $CH^*((GL(2n)/SO(q))_k)\to CH^*((GL(2n)/SO(q))_{\bar{k}})$ is not surjective. $\qed$
\section{$G_2$}
We first record the Chow groups of $G_2$.
\begin{proposition}(Guillot\cite{guillot})
	\[
	CH^*BG_2=\mathbb{Z}[c_2,c_4,c_6,c_7]/(c_2^2=4c_4,c_2c_7=0,2c_7=0)
	\]
\end{proposition} \qed
\\
\\For any field $k$, let $V=k^7$. Let $e_1,\dots, e_7$ denote the standard basis and let $x_1,\dots, x_7$ denote the dual basis of $V^*$. We identify $\bigwedge^3 V^*$ with the space of alternating trilinear forms on $V$. Let $GL(V)$ act on $\bigwedge^3 V^*$ by $(A\cdot f)(x,y,z)=f(Ax,Ay,Az)$.
\begin{lemma}
	\label{3formorbits}
	Let $k$ be a field. The group $GL(V)(\bar{k})$ acts with $9$ orbits on $\bigwedge^3 V^*(\bar{k})$, with representatives $f_1,dots, f_9$ defined over $k$. Let $H_i\subset GL(V)$ denote the stabilizer of $f_i$. Then for all $i<9$, $H_i$ has the mixed Tate property. The orbit of $f_9$ is an open orbit and its stabilizer $H_9$ equals $G_2\times \mu_3$, where $G_2$ denotes the split simple algebraic group of type $G_2$ over $k$.
\end{lemma}
\begin{proposition}
	\label{G2}
	The group $G_2$ has the mixed Tate property over any field $k$.
\end{proposition}
\begin{proof}
	By Lemma \ref{3formorbits} and Corollary \ref{GLnaction}, we conclude that $H_9=G_2\times \mu_3$ has the mixed Tate property. Thus by lemma \ref{mixedtatefunctoriality}a), $G_2$ has the mixed Tate property as well.
\end{proof}
\noindent \textit{Proof of lemma \ref{3formorbits}}: Cohen-Helminck[4] identified the forms $f_i$:
\\$f_1=x_1\wedge x_2\wedge x_3$
\\$f_2=x_1\wedge x_2\wedge x_3+x_1\wedge x_4\wedge x_5$
\\$f_3=x_1\wedge x_2\wedge x_3 + x_4\wedge x_5\wedge x_6$
\\$f_4=x_1\wedge x_6\wedge x_2 + x_2\wedge x_4\wedge x_3 + x_1\wedge x_3\wedge x_5$
\\$f_5=x_1\wedge x_2\wedge x_3 + x_4\wedge x_5\wedge x_6 + x_1\wedge x_4\wedge x_7$
\\$f_6=x_1\wedge x_5\wedge x_2+x_1\wedge x_7\wedge x_4 + x_1\wedge x_6\wedge x_3 + x_2\wedge x_4\wedge x_3$
\\$f_7=x_1\wedge x_4\wedge x_6 + x_1\wedge x_5\wedge x_7 + x_2\wedge x_4\wedge x_5 + x_3\wedge x_6\wedge x_7$
\\$f_8=x_1\wedge x_2\wedge x_3 + x_1\wedge x_4\wedge x_5 + x_1\wedge x_6\wedge x_7$
\\$f_9=x_1\wedge x_2 \wedge x_3+x_4\wedge x_5\wedge x_6+x_1\wedge x_4\wedge x_7 + x_2\wedge x_5\wedge x_7+x_3\wedge x_6\wedge x_7$

Cohen-Helminck proved the above $f_i$ form a set of representatives of the action of $GL(V)(\bar{k})$ on $\bigwedge^3 V^*$ and showed that the stabilizer of $f_9$ is $G_2\times \mu_3$. The forms $f_i$ are defined over $\mathbb{Z}$, hence over $k$. The fact that the orbit of $f_9$ is an open orbit can be seen from the fact that $dim GL(V)/(G_2\times \mu_3) = dim \bigwedge^3 V^*$. 

It remains to prove that the $H_i$ have the mixed Tate property for all $i<9$. Cohen-Helminck calculated the stabilizers $H_i$. We cite the calculations and use them to prove $H_i$ are mixed Tate for $i<9$.
\\\textbf{1.} $H_1=J\rtimes K$, where $J=U'_{3,4}$ and $K=SL(3)\times GL(4)\subset D_{3,4}$. Since $SL(3)$ and $GL(4)$ have the mixed-Tate property, so does $K= SL(3)\times GL(4)$. By lemma \ref{commonexample}, it follows that $H_1$ has the mixed Tate property.
\\
\\\textbf{2.} $H_2=J\rtimes K$, where $J= U'_{1,4,2}$, $K=(1\times Sp(4)\times GL(2))\rtimes G_m\subset D_{1,4,2}$, and $G_m\subset GL(7)$ embeds in the subgroup of diagonal matrices by $\lambda\mapsto (\lambda, \lambda^{-1}, 1, \lambda^{-1}, 1, 1, 1)$.

By Lemma \ref{commonexample}, $N:=J\rtimes K\subset GL(7)$ has the mixed Tate property. Since $H_2/N\cong G_m$, there is a $G_m$-bundle $\pi: GL(7)/N\to GL(7)/H_2$, from which we can conclude that $GL(7)/H_2$ is mixed Tate, so $H_2$ has the mixed Tate property.
\\
\\\textbf{3.} $H_3=J\rtimes K$, where $J=U'_{6,1}$ and $K=((SL(3)\times SL(3))\rtimes \mathbb{Z}/2)\times G_m\subset D_{6,1}$, where $SL(3)\times SL(3)$ embeds diagonally into $GL(6)$ and $\mathbb{Z}/2$ embeds into $GL(6)$ by the permutation $(1 4)(2 5)(3 6)$. Then in fact $(SL(3)\times SL(3))\rtimes \mathbb{Z}/2$ is a wreath product $\mathbb{Z}/2\wr SL(3)$. Since $GL(3)/SL(3)\cong G_m$ is a linear scheme, it follows by Lemma \ref{mixedtatefunctoriality}b' that $K=\mathbb{Z}/2\wr SL(3)$ has the mixed Tate property. So by Lemma \ref{commonexample}, $H_3$ has the mixed Tate property.
\\
\\\textbf{4.} $H_4=J\rtimes K$, where $J\subset U'_{3,3,1}$ consists of those matrices $1+B\in U'_{3,3,1}$ subject to the condition $x_1\wedge Bx_6 \wedge x_2+x_2\wedge Bx_4\wedge x_3 + x_1\wedge x_3\wedge Bx_5$=0;
\\and $K\subset D_{3,3,1}\cong GL(3)\times G_m$, where $GL(3)$ embeds into $D_{3,3,1}$ by $h\mapsto (h,h\cdot det(h)^{-1},1)$ and $G_m$ embeds into $D_{3,3,1}$ by $\lambda\mapsto (1,1,\lambda)$. 

Since $J$ is cut out of $U'_{3,3,1}$ by a linear condition, it is a linear subspace of $U'_{3,3,1}\cong \mathbb{A}^{15}$. Since $K$ has the mixed Tate property, it follows by lemma \ref{commonexample} that $H_4$ has the mixed Tate property. 
\\
\\\textbf{5.} It is more convenient for our purposes to permute the indices, letting $f'_5=x_1\wedge x_3\wedge x_4+x_2\wedge x_5\wedge x_6 +x_1\wedge x_2\wedge x_7$. The form $f'_5$ lies in the same orbit as $f_5$, so we may consider $f'_5$ instead as a representative of the orbit. 

Let $H'_5\subset GL(7)$ be the stabilizer of $f'_5$. Then $H'_5=J\rtimes K$, where $J\subset U'_{2,4,1}$ consists of those matrices $1+B\in U'_{2,4,1}$ subject to the condition $x_1\wedge Bx_3\wedge x_4+x_1\wedge x_3\wedge Bx_4+x_2\wedge Bx_5\wedge x_6+x_2\wedge x_5\wedge Bx_6+e_1\wedge x_2\wedge Bx_7=0$;
\\and $K\cong (GL(2)\times GL(2))\rtimes (\mathbb{Z}/2)$, where $GL(2)\times GL(2)$ embeds into $D_{1,1,2,2,1}\subset D_{2,4,1}$ by $(A,B)\mapsto (det(A)^{-1}, det(B)^{-1},A,B,det(A)det(B))$ and $\mathbb{Z}/2$ embeds into $D_{2,4,1}$ by the permutation $\pi_{(1 2)(3 5)(4 6)}$.

Since $J$ is cut out of $U'_{2,4,1}$ by a linear condition, it is a linear subspace of $U'_{2,4,1}\cong \mathbb{A}^{14}$. The group $K$ is a wreath product $\mathbb{Z}/2\wr GL(2)$. By Lemma \ref{mixedtatefunctoriality}b, $K$ has the mixed Tate property. It follows by Lemma \ref{commonexample} that $H'_5$ has the mixed Tate property. 
\\
\\\textbf{6.} $H_6=J\rtimes K$, where $J\subset U'_{1,3,3}$ consists of matrices $1+B\in U'_{1,3,3}$ subject to the condition $x_1\wedge Bx_5\wedge x_2+x_1\wedge Bx_7\wedge x_4+x_1\wedge Bx_6\wedge x_3=0$, and $K=G_m\times SL(3)\subset D_{1,3,3}$, where $G_m\times SL(3)$ embeds into $D_{1,3,3}$ via $(\lambda,A)\mapsto (\lambda, A, \lambda^{-1}(A^T)^{-1})$. 

Since $J$ is cut out of $U'_{1,3,3}$ by a linear condition, it is a subspace of $U'_{1,3,3}\cong \mathbb{A}^{15}$. Since $K$ has the mixed Tate property, it follows by lemma \ref{commonexample} that $H_6$ has the mixed Tate property.
\\
\\\textbf{7.} Let $W_1=<e_1,e_2,e_3>$ and $W_2=<e_4,e_5,e_6>$. Let $q=x_4x_6+x_5x_7$, a quadratic form on $W_2$.

Then $H_7=J\rtimes K$, where $J\subset U_{3,4}$ consists of the matrices $1+B\subset U_{3,4}$ subject to the condition $Bx_1\wedge x_4\wedge x_6+Bx_1\wedge x_5\wedge x_7+Bx_2\wedge x_4\wedge x_5+ Bx_3\wedge x_6\wedge x_7=0$;
\\and $K=SO(W_2,q)\rtimes G_m\subset D_{3,4}$.

The embedding of $SO(q)$ in $D_{3,4}$ uses the exceptional isomorphism $SO(q)\cong (SL(2)\times SL(2))/(\mathbb{Z}/2)$. Specifically, let $SL(2)\times SL(2)$ act on the space of $2\times 2$ matrices by $(g,h)\cdot x=gxh^{-1}$. This action preserves the determinant, so $SL(2)\times SL(2)$ maps to $SO(q)$ by identifying the tuple $(x_4,x_5,x_6,x_7)$ with the matrix $\left[\begin{array}{cc}
x_6 & x_7
\\-x_5 & x_4
\end{array}\right]$, inducing an isomorphism $\varphi:(SL(2)\times SL(2))/\mu_2\cong SO(q)$. Also, let $SL(2)\times SL(2)$ act on the space of traceless $2\times 2$ matrices by $(g,h)\cdot x = \frac{1}{det(g)^2}gxg^{-1}$. By identifying the tuple $(x_1,x_2,x_3)$ with the traceless matrix $\left[\begin{array}{cc}
x_1 & -x_2
\\x_3 & -x_1
\end{array}\right]$, we obtain a map $\theta: (SL(2)\times SL(2))/\mu_2\to GL(W_1)$. 

Then $SO(W_2,q)$ embeds into $GL(W_1)\times GL(W_2)$ by $A\mapsto (\theta(\varphi^{-1}(A)),A)$. The group $G_m$ embeds into $GL(7)$ diagonally by $\lambda\mapsto (\lambda^{-1},\lambda^{-1},\lambda^{-1},\lambda,1,1,\lambda)$. These two embeddings define the embedding of $SO(W_2,q)\rtimes G_m$ into $D_{3,4}$.

The group $SO(W_2,q)$ has the mixed Tate property since $q$ is a split nondegenerate quadratic form. Thus, $K$ has the mixed Tate property as well. Moreover, $J\subset U_{3,4}$ is a linear subspace, so we conclude by lemma \ref{commonexample} that $H_7$ has the mixed Tate property.
\\
\\\textbf{8.} Let $W=<e_2,\dots,e_7>$ and let $p=x_2\wedge x_3+x_4\wedge x_5+x_6\wedge x_7$, a symplectic form on $W$. 

Then $H_8=J\rtimes K$, where $J= U'_{1,6}$ and $K=(1\times Sp(W,p))\rtimes G_m\subset D_{1,6}$, where $G_m$ embeds into $GL(7)$ diagonally by $\lambda\mapsto (\lambda^{-1},\lambda,1,\lambda,1,\lambda,1)$.

Since $Sp(6)$ has the mixed Tate property, we conclude that $K$ has the mixed Tate property. By lemma \ref{commonexample}, it follows that $H_8$ has the mixed Tate property.
\qed

\section{Quotients by Unipotent groups}
The goal of this section is to prove the following.
\begin{proposition}
	\label{unipotent}
	Let $G$ be an affine group scheme and let $J\subset G$ be a normal split unipotent subgroup scheme of dimension $d$. Suppose $G/J$ acts on a scheme $Z$. Then $M^c([Z/G])\cong M^c([Z/(G/J)])(d)[2d]$. In particular, $M^c(BG)\cong M^c(B(G/J))(d)[2d]$.
\end{proposition}
\noindent Note, the above proposition is not hard to prove in characteristic $0$, but we work over a general field $k$.
\begin{lemma}
	\label{speciallemma}
	Let $B$ be a closed subgroup of $(G_a)^n = \text{Spec } k[x_1,\dots,x_n]$ over $k$. Let $\pi\colon (G_a)^n\to (G_a)^n/B$ be the projection. Let $H((G_a)^n/B)\subset Aut((G_a)^n/B)$ denote the group of automorphisms $\varphi\colon \mathbb{A}^n/B\to \mathbb{A}^n/B$ with the property that for any regular function $f\in k[(G_a)^n/B]$, the degree of $\pi^*f$ equals the degree of $\pi^*\varphi^*f$ in $k[x_1,\dots,x_n]$. Then the group $H((G_a)^n/B)$ is special.
\end{lemma}
\begin{proof}
	Throughout the proof, let $H$ denote the group $H((G_a)^n/B)$.
	
	We first note that $(G_a)^n/B\cong (G_a)^s$, where $s=n-\text{dim}(B)$. Indeed, Rosenlicht proved\cite{ros} that $B\subset (G_a)^{n}$ is cut out by $s$ additive functions $y_1,\dots,y_s\colon (G_a)^{n}\to G_a$ which are algebraically independent. Then $y:=y_1\times \cdots \times y_s\colon (G_a)^{n}\to (G_a)^s$ is a dominant homomorphism with kernel $B$, showing that $(G_a)^n/B\cong (G_a)^s$. 
	
	Let $I\subset k[x_1,\dots, x_n]$ denote the ideal cutting out $B$ from $(G_a)^n$. Let $p=$char $k$. Then each $y_i\in I\subset k[x_1,\dots, x_n]$, being additive, must be a $p$-polynomial of the $x_i$, that is a linear combination of the monomials $x_i^{p^r}$ (If $p=0$, a $p$-polynomial means a linear combination of the monomials $x_i$). 
	
	Let $p^m$ be the highest degree out of all the polynomials $y_i\in I\subset k[x_1,\dots,x_n]$, and let $V\subset I$ be the vector space of all $p$-polynomials in $I$ whose degree is at most $p^m$. Note that every polynomial $y\in I$ is of the form $\pi^*(z)$ for a unique regular function $z\in k[\mathbb{A}^n/B]$. Then $V$ is a faithful representation of $H$, where an automorphism $\varphi\in H$ acts by mapping each $y\in V$ to $\pi^*(\varphi^*(z))$, for $z\in k[\mathbb{A}^n/B]$ such that $\pi^*(z)=y$.
	
	An automorphism in $H$ is determined by its action on the $y_i$, since any additive polynomial in $I$ is a $p$-polynomial of the $y_i$. Thus $H\subset GL(V)$ consists precisely of $h\in GL(V)$ preserving degree and such that $h\cdot (v^p)=(h\cdot v)^p$ whenever $v,v^p\in V$.
	\\
	\\We can choose the $y_i$ such that the top-degree terms of all the $y_i^{p^{r}}$, ranging over all $i$ and all $r\in \mathbb{N}$ are linearly independent. Indeed, suppose we have a nontrivial linear relation between the top-degree terms. We can assume that the polynomials $y_i^{p^{r}}$ involved in this relation all have the same degree $p^d$. This implies that there exists a nontrivial linear combination $\sum c_i y_i^{p^{r_i}}$ of degree less than $p^d$. Let $r$ be the smallest of the $r_i$ such that $c_i\ne 0$. Note the top-degree terms of $y_i^{p^{r_i}}$ are linearly independent if and only if the top-degree terms of $y_i^{p^{r_i-r}}$ are linearly independent. Thus we may reduce to the case of $r=0$.
	
	It follows that, for some $i$, a polynomial of the form $y'_i=y_i-\sum_{j\ne i} c_j y_j^{p^{r_j}}$ has smaller degree than $y_i$. We may replace $y_i$ with $y'_i$ and still retain an additive generating set of $I$. Since this process reduces the degrees of the $y_i$, it must eventually terminate, and on terminating the $y_i$ must have the desired property.
\end{proof}
\noindent We now introduce two flags on $V$. On the one hand, for each $0\le i\le m$, let $F_i\subset V$ be the subspace of all polynomials $V$ of degree at most $p^i$. Let $F'_i\subset V$ be the subspace spanned by all the $(y_j)^{p^r}$ of degree $p^i$. Then, by our choice of the $y_j$ to have the property that the top-degree terms of $(y_j)^{p^r}$ are linearly independent, it follows that $F_i=\oplus_{j\le i}F'_j$. Since automorphisms in $H$ preserve degree, $H$ preserves each $F_i$, though $H$ may not preserve the $F'_i$. Conversely, if $g\in GL(V)$ preserves the $F_i$, then $g$ must also preserve degree, since for each $v\in F'_i$, the component of $gv$ in $F'_i$ must be nonzero by nonsingularity, so $g$ preserves the degree of vectors in each $F'_i$, hence in all of $V$.

On the other hand, for each $0\le i\le m$, let $W'_i\subset V$ be spanned by all the $(y_j)^{p^r}\in V$ such that $y_j$ has degree exactly $p^i$. Let $W_i=\oplus_{j\le i} W'_i$. Then $H$ preserves each $W_i$ since the action of $H$ on $V$ preserves degree and commutes with taking $p$th powers. However, $H$ may not preserve the $W'_i$.

Let $T\subset GL(V)$ be the subgroup preserving the flag $W_i$. Let $D\subset T$ be the subgroup preserving each of the $W'_i$. Let $\pi\colon T\to D$ be the natural projection, and let $U$ be its kernel. We claim that $\pi$ maps $H$ to $H$.

Indeed, let $h\in H$, $v\in W_i$. Then $h\cdot v=\pi(h)\cdot v+v'$, where $\pi(h)\cdot v\in W'_i$, $v'\in W_{i-1}$. By our choice of the $y_i$, it follows that $\pi(h)\cdot v$ has degree no larger than the degree of $h\cdot v$, which equals the degree of $v$. Thus $\pi(h)$ preserves each $F_i$, so by an earlier observation, $\pi(h)$ preserves degree. Furthermore, since $(W'_i)^p=\subset W'_{i}$ for each $i$, it follows that $\pi(h)\cdot v^p=(\pi(h)\cdot v)^p$. So $\pi(h)\in H$, as desired.

Thus, $H$ splits as the semidirect product of $H\cap D$ and $H\cap U$. Let $W''_i\subset W'_i$ be the subspace spanned by all the $y_j$ of degree $p^i$. An automorphism of $W'_i$ commuting with taking $p$th powers precisely corresponds to an automorphism of $W''_i$. Hence, $H\cap D\cong \Pi_i GL(W''_i)$.

Now, consider the standard normal series $U_i$ of $U$, where $U_i$ consists of all $h\in U$ such that, for each $j>i$ and $v\in W'_j$, the projection of $h(v)$ onto $W_{j-i-1}$ is $0$. The groups $U_i\cap H$ form a normal series of $U\cap H$. Similarly to above, we have an isomorphism $(U_{i+1}\cap H)/(U_i\cap H)\cong \Pi_j \text{Hom}(W''_j,F_j\cap W'_{j-i})$ as commutative groups. We conclude that $U\cap H$ is a split unipotent group. 

Thus $H$ is an extension of copies of $GL(n_i)$ for various $n_i$ and copies of $G_a$. Since extensions of special groups are special, and $GL(n_i)$, $G_a$ are special, it follows that $H$ is special. $\qed$
\begin{lemma}
	\label{flag}
	Let $V$ be a vector space over $k$. Let $G\subset GL(V)$ be a subgroup scheme and let $J\subset G$ be a normal unipotent subgroup scheme. Then there exists a flag $0=V_0\subset V_1\subset V_2\subset \cdots V_k=V$ which is preserved by $G$ and such that $J$ acts as the identity on each quotient space $V_j/V_{j-1}$ for $j=1,\cdots, k$.
\end{lemma}
\begin{proof}
	By induction, it suffices to find a subspace $V_1\subset V$ preserved by $G$, on which $J$ acts trivially. We use the fact that $J$ is unipotent to find $v\in V$ such that $jv=v$ for all $j\in J$. Let $V_1$ be the subspace spanned by all vectors $gv$, for all $g\in G$. The group $G$ clearly preserves $V_1$. To show $J$ acts as the identity on $V_1$, we use the fact that, for any $j\in J$, $g\in G$, we have $g^{-1}jg\in J$ since $J$ is normal in $G$. In particular, $g^{-1}jgv=v$, so $jgv=gv$, so $j$ fixes all vectors $gv$, as desired.
\end{proof}
\noindent In the context of the above lemma, let $n_i=\text{dim}(V_i/V_{i-1})$. It follows that there is a choice of basis $e_i$ for $V$ such that $G$ is contained in the group $T_{n_1,\dots,n_k}$ fixing the flag $\langle e_1,\dots,e_{n_1}\rangle\subset \cdots \subset \langle e_1,\dots,e_{n_1+\cdots + n_i}\rangle \subset \cdots \subset V$, and $J$ is contained in the subgroup $U_{n_1,\dots,n_k}$ of $T_{n_1,\dots, n_k}$ acting trivially on each quotient space $\langle e_1,\dots,e_{n_{1}+\cdots+n_i+n_{i+1}}\rangle / \langle e_1,\dots, e_{n_1+\cdots+n_i} \rangle$.
\begin{lemma}
	\label{fullunipotent}
	Let $n=n_1+\cdots + n_k$. Let $U:=U_{n_1,\dots,n_k}, T:=T_{n_1,\dots,n_k}\subset GL(n)$. Suppose $G\subset T$ and $J\subset G\cap U$ is normal in $G$. Let $X$ be a scheme on which $T$ acts freely and let $Y$ be a scheme on which $G/J$ acts freely. Then the $U/J$-fibration $(X/J\times Y)/(G/J)\to (X/U\times Y)/(G/J)$ induces by flat pullback an isomorphism $M^c((X/U\times Y)/(G/J))(r)[2r]\cong M^c((X/J\times Y)/(G/J))$, where $r=\text{dim}(U/J)$.
\end{lemma}
\begin{proof}
	Note that $G$ normalizes $U$, so $G/J$ indeed acts on $X/U$. 
	
	Let $1=Z_0\subset Z_1\subset \cdots \subset Z_m=U$ be a central normal series for $U$, whose successive quotients are isomorphic to $(G_a)^{r_i}$ for various $r_i$. We can furthermore ensure that each $Z_i$ is normal in $T$, by letting $Z_i$ be cut out of $Z_{i+1}$ by setting all entries in a given $n_j\times n_k$-block equal to $0$.
	
	Let $E_i=(X/Z_i)/(J/(Z_i\cap J))$. Note we have a sequence of maps 
	\\$X/J=E_0\to E_1\to \cdots \to E_{m-1}\to E_m=X/U$. Each map $E_i\to E_{i+1}$ is a principal bundle with fiber $A_i/B_i$, where $A_i:=(Z_{i+1}/Z_i),B_i:=(J/(J\cap Z_i))$.
	
	Since $G$ normalizes each of the $Z_i$, we have a sequence of $A_i/B_i$-fibrations:
	\\$(X/J\times Y)/(G/J)\to (E_1\times Y)/(G/J)\to \cdots \to (E_{m-1}\times Y)/(G/J)\to$
	\\$(X/U\times Y)/(G/J)$.
	
	Each of the above fibrations has a flat trivialization:
	
	\[
	\xymatrix{
		E_i\times Y\times A_i/B_i \ar[r] \ar[d] & E_i\times Y \ar[d]
		\\
		E_i\times Y \ar[r] \ar[d]^{p_i} & E_{i+1}\times Y \ar[d]^{p_{i+1}}
		\\
		(E_i\times Y)/(G/J) \ar[r]^{\pi_i} & (E_{i+1}\times Y)/(G/J)
	} 
	\]
	Consider a base point $b\in (E_{i+1}\times Y)/(G/J)$. Let $F_b= \pi_i^{-1}(b)$. Each point $(e,y)\in E_{i}\times Y$ over $b$ gives a trivialization $\psi_{e,y}\colon A_i/B_i\cong F_b\colon a\mapsto p_i(ae,y)$, for $a\in A_i/B_i$. Now, all other points in $E_i\times Y$ over $b$ are of the form $(a'ge,gy)$, where $g\in G/J, a'\in A_i/B_i$. We have $p_i(ae,y)=p_i(gae,gy)=p_i(((gag^{-1})a'^{-1}\cdot a'ge,gy))$. Thus the transition function $\psi_{a'ge,gy}^{-1}\circ \psi_{e,y}$ is $a\mapsto (gag^{-1})a'^{-1}$.
	
	Recall $A_i=(Z_{i+1}/Z_i)\cong (G_a)^{r_i}$ corresponds to a certain $n_j\times n_k$-block of matrices in $U$. There is a natural vector space structure on $A$ with respect to which conjugation by an element of $G$ is linear, since it corresponds to multiplication by fixed matrices, and translation by an element of $A_i$ is affine-linear since it corresponds to addition of a fixed matrix. Thus the transition functions on $A_i/B_i$ descend from affine-linear transformations on $A_i$. In particular, each transition function $\varphi\colon A_i/B_i\to A_i/B_i$ has the property that, for each $f\in k[A_i/B_i]$, the degree of $q^*f$ equals the degree of $q^*\varphi^*f$. Thus the transition functions belong to the group $H(A_i/B_i)\subset \text{Aut}(A_i/B_i)$ described in Lemma 1. Since $H(A_i/B_i)$ is special, it follows that each fibration $\pi_i$ is Zariski-locally trivial with fiber $A_i/B_i$ isomorphic to affine space. Thus pullback by $\pi_i$ induces an isomorphism of compactly supported motives (with a shift), so the composite pullback is an isomorphism as well.
\end{proof}
\noindent Recall that given an affine group scheme $G$ we say that $(V_i-S_i)/G$ is a resolution of $BG$ by algebraic spaces if the $V_i$ are representations of $G$ of dimensions $n_i$, $S_i\subset V_i$ are $G$-invariant closed subschemes outside which $G$ acts freely, such that the codimension of $S_i$ in $V_i$ approaches infinity and such that we are given linear surjections $f_i\colon V_{i+1}\to V_i$ with $S_{i+1}\subset f_i^{-1}(S_i)$. We think of such as a resolution as an object of the category $\mathcal{A}$ introduced in Section 3.
\\
\\
\textit{Proof of Proposition \ref{unipotent}}: Let $G\subset GL(n)$ be a subgroup scheme and let $J\subset G$ be a normal unipotent subgroup scheme. By Lemma 2, we may assume after change of basis that $J\subset U:=U_{n_1,\dots,n_k}$ and $G\subset T:=T_{n_1,\dots,n_k}$.
\\
\\Let $(V_i-S_i)/T$ be a resolution of $BT$ by algebraic spaces. Then consider the corresponding resolution $X_i:=(V_i-S_i)/J$ of $BJ$ by algebraic spaces. Then $G/J$ acts freely on $X_i$. Let $n$ be the relative dimension of $X_i/(G/J)\to B(G/J)$ in $\mathcal{A}$.  We claim that $(Z\times X_i)/(G/J)$ satisfies the hypotheses of Lemma \ref{generalresolution} with respect to the quotient stack $[Z/(G/J)]$. Once we show this, the proposition will follow, since by Lemma \ref{generalresolution}, $M^c([Z/(G/J)])$ will be isomorphic to the homotopy limit of $LM^c((X_i\times Z)/(G/J))(-n)[-2n]$. But $(X_i\times Z)/(G/J)\cong ((V_i-S_i)\times Z)/G$, so the above homotopy limit is isomorphic to $M^c([Z/G])(-d)[-2d]$.
\\
\\Let $Y$ be an algebraic space on which $G/J$ acts. As in Lemma \ref{generalresolution}, let $H(Y)$ denote the (canonical) homotopy limit $LM^c(((V_i-S_i)/J\times Y)/(G/J))(-n)[-2n]$. Concretely, the hypothesis of Lemma \ref{generalresolution} amounts to the statement that the pullback $M^c(Y/(G/J))\to H(Y)$ be an isomorphism for all $Y$. Let $\tilde{H}(Y)$ denote the similar homotopy limit $LM^c(((V_i-S_i)/U\times Y)/(G/J))(-n+r)[-n+2r]$, where $r=\text{dim}U/J$. By Lemma \ref{fullunipotent}, flat pullback by $((V_i-S_i)/J\times Y)/(G/J)\to ((V_i-S_i)/U\times Y)/(G/J)$ induces an isomorphism of compactly supported motives with dimension shifts, so the homotopy limit of these isomorphisms yields a canonical isomorphism $\tilde{H}(Y)\cong H(Y)$. Thus we can reduce the question to showing that $\tilde{H}(Y)\cong M^c(Y/(G/J))$.
\\
\\
First consider the principal $T/U$-bundle $p\colon (T/U\times Y)/(G/J)\to Y/(G/J)$. Since $T/U\cong \Pi GL(n_i)$ is special, $p$ is Zariski-locally trivial. So let $O\subset Y$ be a $(G/J)$-invariant open set such that $(T/U\times O)/(G/J)\to O/(G/J)$ is a trivial $T/U$-bundle, with trivializaton $\varphi\colon (T/U\times O)/(G/J)\cong T/U\times O/(G/J)$. We claim that each fibration $\pi_i\colon ((V_i-S_i)/U\times Y)/(G/J)\to Y/(G/J)$ is also trivial over $O$. 

Let $q\colon (T/U\times O)/(G/J)\to T/U$ be the composition of $\varphi$ with projection onto the first component. Since $T$ acts on $V$, $T/U$ acts on $(V_i-S_i)/U$. We then define a map $\psi\colon (V_i-S_i)/U\times O\to (V_i-S_i)/U\times O/(G/J)$ by $(v,y)\mapsto (q(\overline{1,y})v,\overline{y})$, for $v\in (V_i-S_i)/U, y\in O$. Let $g\in G/J$ and for simplicity let $g$ also denote the image of $g$ in $T/U$. Note that $\psi(gv,gy)=(q(\overline{1,gy})gv,\overline{gy}))=(q(\overline{g^{-1},y})gv,\overline{y})=(q(\overline{1,y})g^{-1}gv,\overline{y}))=\psi(v,y)$, where we have used the $T/U$-equivariance of the trivialization $\varphi$. Since $\psi$ is constant on $G/J$-orbits, it descends to give the desired trivialization of the bundle $\pi_i$ over $O$.
\\
\\Now, let $Z$ be the closed complement of $O$ in $Y$. We have a map of triangles:
\[
\xymatrix{
	M^c(Z/(G/J)) \ar[r] \ar[d]^{\alpha} & M^c(Y/(G/J)) \ar[r] \ar[d]^{\beta} & M^c(O/(G/J)) \ar[d]^{\gamma}
	\\
	\tilde{H}(Z) \ar[r] & \tilde{H}(Y) \ar[r] & \tilde{H}(O)
} 
\]
where the vertical maps are homotopy limits of flat pullback maps.

The map $\gamma$ is an isomorphism. Indeed, we know that each of the schemes $((V_i-S_i)/U\times O)/(G/J)$ is isomorphic to $(V_i-S_i)/U \times O/(G/J)$. The projection morphism $(V_i-S_i)\times O/(G/J)\to (V_i-S_i)/U \times O/(G/J)$ is a principal $U$-bundle, hence it induces an isomorphism of compactly supported motives with shifts. Thus $\tilde{H}(O)$ is isomorphic to the homotopy limit of the motives $M^c((V_i-S_i)\times O/(G/J))(-l_i)[-2l_i]$, where $l_i=\text{dim}(V_i)$. But by the localization triangles, and since the codimension of $S_i$ in $V_i$ approaches infinity, this homotopy limit is $M^c(O/(G/J))$.

By induction on the size of an open cover for the trivialization, $\alpha$ is also an isomorphism. It follows that $\beta$ is an isomorphism, as desired. $\qed$
\bibliography {main}    % bibliography references
\bibliographystyle {alpha}

\end {document}